\documentclass[a4paper]{article}
\usepackage{amsmath}
\usepackage{amsthm}
\usepackage{amssymb}
\usepackage{graphicx}
\usepackage{authblk}
\usepackage{hyperref}
\usepackage{multirow}
\usepackage{cite}
\usepackage{comment}

\newtheorem{theorem}{Theorem}[section]
\newtheorem{proposition}[theorem]{Proposition}
\newtheorem{lemma}[theorem]{Lemma}

\newtheorem{definition}{Definition}
\newtheorem{remark}{Remark}

\newenvironment{vf}{\left\{\begin{array}{rcl}}{\end{array}\right.}

\DeclareMathOperator{\cycl}{Cycl}

\newcommand{\cdim}{\mathop{\mathrm{co\textrm{-}dim}}}
\newcommand{\includegraph}[2][]{\ifnum\pdfoutput=0\includegraphics[#1]{#2.eps}\else\includegraphics[#1]{#2.pdf}\fi}

\author[1]{Peter De Maesschalck}
\author[2]{Renato Huzak}
\author[3]{Ansfried Janssens\footnote{Corresponding author, {\tt ansfried.janssens@uhasselt.be}}}
\author[4]{Goran Radunovi\'{c}}

\affil[1,2,3]{Hasselt University, Campus Diepenbeek, Agoralaan Gebouw D, 3590 Diepenbeek, Belgium}
\affil[4]{University of Zagreb, Faculty of Science, Horvatovac 102a, 10000 Zagreb, Croatia}

\title{Fractal codimension of nilpotent contact points in two-dimensional slow-fast systems}

\date{}

\begin{document}
\maketitle

\begin{abstract}
 In this paper we introduce the notion of fractal codimension of a nilpotent contact point $p$, for $\lambda=\lambda_0$, in smooth planar slow–fast systems $X_{\epsilon,\lambda}$ when the contact order $n_{\lambda_0}(p)$ of $p$ is even, the singularity order $s_{\lambda_0}(p)$ of $p$ is odd and $p$ has finite slow divergence, i.e., $s_{\lambda_0}(p)\le 2(n_{\lambda_0}(p)-1)$. The fractal codimension of $p$ is a generalization of the ``traditional'' codimension of a slow-fast Hopf point of Li\'{e}nard type, introduced in (Dumortier and Roussarie (2009), \cite{DRbirth}), and it is intrinsically defined, i.e., it can be directly computed without the need to first bring the system into its normal form. The intrinsic nature of the notion of fractal codimension stems from the Minkowski dimension of fractal sequences of points, defined near $p$ using the so-called entry-exit relation, and slow divergence integral. We apply our method to a slow-fast Hopf point and read its degeneracy (i.e., the first nonzero Lyapunov quantity) as well as the number of limit cycles near such a Hopf point directly from its fractal codimension. We demonstrate our results numerically on some interesting examples by using a simple formula for computation of the fractal codimension.
\end{abstract}

\tableofcontents

\section{Introduction}\label{Section-Introduction}
In the study of ($C^\infty$-)smooth planar slow-fast vector fields one often encounters a slow-fast Hopf point (or a singular Hopf
point) from which limit cycles may bifurcate. Analysis of limit cycles near such a Hopf point is called the birth of canards and is usually done in the following normal form near the origin $(x,y)=(0,0)$:
 \begin{equation}
\label{model-Hopf}
    \begin{vf}
        \dot{x} &=& y-x^2+x^3h_1(x,\lambda)  \\
        \dot{y} &=&\epsilon\big(b(\lambda)-x+x^2h_2(x,\epsilon,\lambda)+yh_3(x,y,\epsilon,\lambda)\big),
    \end{vf}
\end{equation}
where $\epsilon\ge 0$ is the singular perturbation parameter which is kept small, $\lambda\sim\lambda_0\in\mathbb R^l$, $b(\lambda_0)=0$ and $b$, $h_1$, $h_2$ and $h_3$ are smooth functions. For more details see Section 6.1 in \cite{DDR-book-SF}. The number of limit cycles of \eqref{model-Hopf} in a small $(\epsilon,\lambda)$-uniform neighborhood of the slow-fast Hopf point $(x,y)=(0,0)$, with $(\epsilon,\lambda)\sim (0,\lambda_0)$, typically depends on codimension of the slow–fast Hopf point (the higher the codimension, the more limit cycles can be born). A natural question that arises is how to define the notion of codimension of the Hopf point in \eqref{model-Hopf}. Assume that $h_3=0$ (thus, \eqref{model-Hopf} is a Li\'{e}nard system). The birth of canard cycles in this Li\'{e}nard setting has been studied in \cite{DRbirth}. Following \cite{DRbirth}, it is possible to eliminate $h_2$ by means of a smooth equivalence. Now, in classical Li\'{e}nard setting ($h_2=h_3=0$) the notion of codimension is defined as follows (see \cite{DRbirth}): the singular Hopf point has codimension $j+1\ge 1$ if 
\[h_1(x,\lambda_0)+h_1(-x,\lambda_0)=\alpha x^{2j}+O(x^{2j+2}), \ \ \alpha\ne 0.\]
If the slow-fast Hopf point has finite codimension $j+1$, then the maximum number of limit cycles produced by the slow-fast Hopf point is finite (see \cite{DRbirth}) and, moreover, bounded by $j+1$ (see \cite{Li-conjecture-solved} and Remark \ref{remark-conjecture} in Section \ref{section-statement}). The codimension one case has been studied in \cite{1996,KS}.
\smallskip

If we deal with analytic families in \eqref{model-Hopf}, then we can eliminate $h_3$ reducing the study of the birth of canards to Li\'{e}nard form (see \cite{RHAnalytic}). But, as explained above, even 
in the Li\'{e}nard setting determination of the codimension of slow-fast Hopf points depends on normal form of classical Li\'{e}nard type. \textit{Since putting systems into normal form is often very challenging from computational point of view, in this paper we present an intrinsic fractal approach for a direct determination of so-called fractal codimension of slow-fast Hopf points or more general contact points (see Definition \ref{def-intrin-codimen} in Section \ref{section-statement}). There is no need to use normal forms. } In classical Li\'{e}nard setting, this fractal codimension coincides with the codimension introduced in \cite{DRbirth} (see Lemma \ref{lemma-trad-frac} in Section \ref{section-statement}). 
\smallskip

Following \cite{DDR-book-SF}, there are important notions that can be intrinsically defined inside two-dimensional slow-fast systems, for example fast foliation, slow-fast Hopf point, contact point, slow divergence integral, etc. This means that their definition is coordinate-free (see Section \ref{section-statement}). We also refer to \cite{Martin} for a coordinate-free approach. Using the intrinsic nature of the notion of a slow-fast Hopf point, in \cite{CriticalityHopf} it has been shown that the criticality of a slow-fast Hopf bifurcation can be intrinsically determined using a simple formula. The paper covers the case when the first Lyapunov quantity is different from zero. More degenerate slow-fast Hopf points are more difficult to treat using the coordinate-free approach introduced in \cite{CriticalityHopf} because formulas become quite long as we increase the codimension (in the sense of \cite{DRbirth}). \textit{The main advantage of our fractal (coordinate-free) approach is that the fractal codimension can be determined using a single closed formula given in Definition \ref{def-intrin-codimen}.} Our result covers degenerate (or Bautin) slow-fast Hopf bifurcations of any fractal codimension, including $\infty$ for analytic systems (Theorem \ref{theorem-Lienard} in Section \ref{section-statement}).
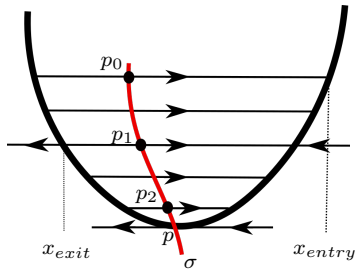
\begin{figure}[htb]
	\begin{center}
		\includegraphics[width=4.5cm,height=3.3cm]{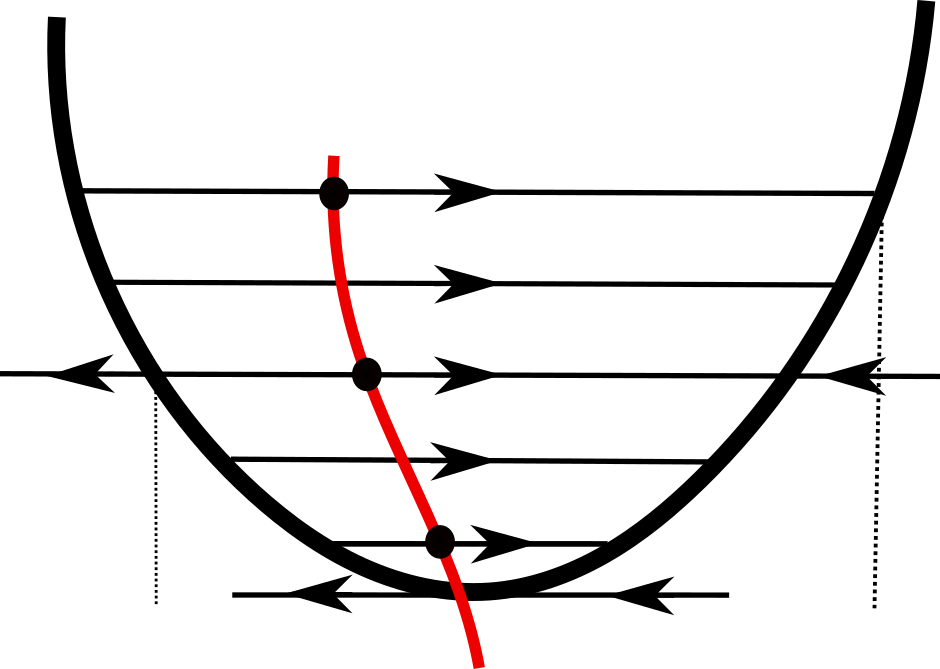}
		{\footnotesize
        \put(-115,-0){$x_{exit}$}\put(-21,-0){$x_{entry}$}\put(-62,-5){$\sigma$}
        \put(-70,5){$p$}
        \put(-93,71){$p_0$}
        \put(-89,43){$p_1$}
        \put(-80,21){$p_2$}
}
         \end{center}
	\caption{Dynamics of \eqref{model-Hopf} near the origin, for $\epsilon=0$ and $\lambda=\lambda_0$, and a section $\sigma$, in a small neighborhood of $p$, transversally cutting the critical curve at $p$. $\mathcal S=\{p_k \ | \ k\ge 0\}$ is a fractal sequence.}
	\label{fig-chirp1}
\end{figure}
\smallskip

Let us now explain the fractal approach mentioned above in more detail. We define fractal sequences of points and associated geometric chirps \cite{Goran} using the so-called entry-exit relation and the ``tunnel'' behaviour near a slow-fast Hopf point $p=(0,0)$ (see \cite{Benoit,DM-entryexit}). When $\epsilon=0$ and $\lambda=\lambda_0$ in \eqref{model-Hopf}, the point $p$ separates two branches of the critical curve $y=x^2(1+O(x))$, one normally attracting ($x>0$), the other normally repelling ($x<0$). We focus on a portion of the fast foliation, given by horizontal lines, located in between the two branches, and choose a regular section $\sigma$ transverse to the portion and containing $p$  (see Fig. \ref{fig-chirp1}). Roughly speaking, a fractal sequence is the union of countably infinitely many points on section $\sigma$ accumulating at $p$, where one jumps from one point to the next one using the entry-exit relation:
\[  \int_{x_{entry}}^{x_{exit}} 4x(1+O(x))dx=0, \ \ x_{entry}>0,\ x_{exit}<0.      \]
We call this integral the slow divergence integral along the critical curve between $x_{entry}$ and $x_{exit}$ (see \cite{DDR-book-SF} or Section \ref{section-statement} for more details). We point out that the integrand changes sign as $x$ varies through the origin and, as a simple consequence of this, for each $x_{entry}$ small enough we can find a unique $x_{exit}$ such that the entry-exit relation holds (in the literature, this is sometimes called the ``tunnel'' phenomenon if there exists the so-called breaking parameter). More precisely, we define a fractal sequence $\mathcal S=\{p_k \ | \ k\ge 0\}$ as follows: we start with a $p_0\in\sigma$ and denote by $x=x_{entry}$ the $x$-component of the $\omega$-limit point of $p_0$. Then we find $p_1\in\sigma$ such that the $x-$component of the $\alpha$-limit point of $p_1$ is $x=x_{exit}$, using the entry-exit relation. We denote the $x$-component of the $\omega$-limit point of $p_1$ again by $x=x_{entry}$ and define $p_2\in \sigma$ using the entry-exit relation, etc. It is clear that near the slow-fast Hopf point we can define different fractal sequences--by simply changing the level of the starting point $p_0$ in the iterative definition of fractal sequence or by changing $\sigma$. We refer to Section \ref{section-statement} for precise assumptions on \eqref{model-Hopf} under which such fractal sequences are well-defined (i.e. converge monotonically to $p$).  For a fixed smooth section $\sigma$ transversally cutting the critical curve at $p$ and for a fixed fractal sequence $\mathcal S$ on $\sigma$ we denote by $\mathcal C_{\mathcal S}$ the union of fast orbits passing through the points of $\mathcal S$, bounded by the critical curve, and call it the associated chirp.
\smallskip

We want to ``measure the density'' of $\mathcal{S}$ and $\mathcal{C}_{\mathcal S}$ using the notion of Minkowski dimension (equivalent to the box dimension). See e.g. \cite{Falconer,Goran} and Section \ref{section-Minkowski-def}. \textit{The bigger the Minkowski dimension of $\mathcal{S}$ and $\mathcal{C}_{\mathcal S}$, the higher the density of $\mathcal{S}$ and $\mathcal{C}_{\mathcal S}$, and the fractal codimension of the slow-fast Hopf point increases} (see Section \ref{section-statement}). We will show that the Minkowski dimension of $\mathcal{S}$ and $\mathcal{C}_{\mathcal S}$ is well-defined, i.e. independent of the chosen $\sigma$ and the starting point $p_0\in\sigma$ (see Theorem \ref{theorem-Main} in Section \ref{section-statement}).  
\smallskip

The contact at the slow-fast Hopf point in \eqref{model-Hopf} between the critical curve and the
leaf through $p=(0,0)$ of the fast foliation, at level $(\epsilon,\lambda)=(0,\lambda_0)$, is of order $2$ and singularity order at $p$ is $1$ (the singularity order has more involved definition, see Section \ref{section-statement}). We will work in a more general framework in which smooth planar nilpotent contact points of arbitrary even contact order and odd singularity order are allowed. There is an additional restriction: such contact points need to have finite slow divergence integral. For more details see Assumptions A, B, C and C' in Section \ref{section-statement}. Then we can use an entry-exit relation and (intrinsically) define $\mathcal S$ and $\mathcal C_{\mathcal S}$ in a similar fashion. We will show that the Minkowski dimension of $\mathcal S$ exists and can
take only a discrete set of values given in Theorem \ref{theorem-1} or Theorem \ref{theorem-Main}. (We get a partial result for the Minkowski dimension of $\mathcal C_{\mathcal S}$.) Roughly speaking, with the above assumptions Theorem \ref{theorem-Main} states the following:
\begin{itemize}
    \item The Minkowski dimension $\dim_B\mathcal S$ of $\mathcal S$ exists and \begin{equation}\label{dimS-j}
  \dim_B\mathcal S\in \left\{\frac{2j+1}{n+2j+1} \ | \ j\in\mathbb N_0 \right\}\cup \{1\},  
\end{equation}
where a nilpotent contact point has  even contact order $n\ge 2$. The Minkowski dimension of $\mathcal S$ is a coordinate-free notion and it does not depend on the choice of the section $\sigma$ and the first element $p_0$ of $\mathcal S$.
\end{itemize}

The coordinate-free nature of the Minkowski dimension and its bijective correspondence to $j$ appearing in \eqref{dimS-j} enable us to define the notion of fractal codimension of the nilpotent contact point $p$ (Definition \ref{def-intrin-codimen}):
\begin{itemize}
    \item  If $\dim_B\mathcal S<1$, we say that a nilpotent contact point has finite fractal codimension, $$\cdim(p):=j+1\ge 1$$
with
\[j=\frac{(n+1)\dim_B\mathcal S-1}{2(1-\dim_B\mathcal S)}\in\mathbb{N}_0.\]
When $\dim_B\mathcal S=1$, then we say that the fractal codimension is infinite.
\end{itemize}

\smallskip

Furthermore, in Theorem \ref{theorem-Lienard} we show that the number of limit cycles produced by an analytic slow-fast Hopf point is bounded from above by the fractal codimension, i.e. it cannot exceed $j+1$. Moreover, the same is true in the smooth setting if, in addition, a slow-fast Hopf point allows reduction to a classical Li\'{e}nard normal form for smooth equivalence.

\smallskip 

A similar fractal approach has been used in \cite{BoxRenato,BoxVlatko,BoxDomagoj} to study the Minkowski dimension of so-called balanced canard cycles in planar slow-fast families. We refer to \cite{DDR-book-SF} for more details about such canard cycles which are of size $O(1)$ in the phase space. A fractal analysis of planar nilpotent singularities (without parameters) using unit time map can be found in \cite{BoxVesna}. We also refer to \cite{ZupZub} where the Minkowski dimension of spiral trajectories of some planar vector fields has been studied.
\smallskip

In this paper we deal with two-dimensional slow-fast systems defined on a smooth surface (e.g., $\mathbb R^2$). It is natural to consider smooth slow-fast systems in this setting. On one hand, we can directly use normal forms for smooth equivalence near contact points from \cite{DDR-book-SF}. On the other hand, we will show that finite degrees of smoothness are not sufficient when the Minkowski dimension of $\mathcal S$ tends to $1$ (i.e., when the fractal codimension tends to infinity).

\smallskip

We refer to
\cite[Section 16.2]{DDR-book-SF} for another definition of codimension of contact points, different from our ``fractal'' definition.

\smallskip

The paper is organized as follows.
In Section \ref{section-Minkowski-def} we recall some basic properties of the Minkowski dimension and in Section \ref{subsection-riemann} we show that it is independent of the choice of the Riemannian metric on a smooth surface. In Section \ref{section-statement} we state our main results. We prove the main results in Section \ref{section-proofs}.
Finally, section \ref{section-applications} is devoted to applications.

\section{The Minkowski dimension and its properties}
\label{section-Minkowski-def}

The fractal tool that we will use to determine the first nonzero Lyapunov or saddle quantity is the notion of Minkowski dimension (always equal to the well known box dimension).\footnote{Note that in our context the more popular Hausdorff dimension is not applicable since due to its \emph{countable stability property} it would be trivial on the sets generated by the dynamics of the slow-fast system.}
Let us briefly recall the definitions and some of its main properties.

Let $A$ be a nonempty compact subset of the Euclidean space $\mathbb{R}^N$.
By $A_{\varepsilon}$ one denotes its \emph{$\varepsilon$-neighborhood} in $\mathbb{R}^N$, i.e.,
\begin{equation}
	\label{eps-neigh}
	A_{\varepsilon}:=\{x\in\mathbb{R}^N|d(x,A)<\varepsilon\},
\end{equation}
where $d(x,A)$ is the Euclidean distance from the point $x$ to the set $A$.
Then, for $r\geq0$, one defines the \emph{$r$-dimensional Minkowski content of $A$} as the following limit (if it exists):
\begin{equation}
	\label{r-Mink}
	\mathcal{M}^r(A):=\lim_{\varepsilon\to0^+}\frac{|A_{\varepsilon}|}{\varepsilon^{N-r}},
\end{equation}
where $|A_{\varepsilon}|$ denotes the $N$-dimensional Lebesgue measure of $A_{\varepsilon}$, i.e., the volume of the $\varepsilon$-parallel body of $A$.
Roughly speaking, one wants to extract the leading term in the asymptotics of the so-called \emph{tube function of the set $A$}: $\varepsilon\mapsto|A_{\varepsilon}|$, i.e., the leading term of its Steiner-like formula.
For ``nice'' sets this is possible and then the critical value of $r=D$ for which $\mathcal{M}^D(A)$ is positive and finite is called the \emph{Minkowski dimension} of the set $A$ and denoted by $D=\dim_BA$ and, moreover, the set is said to be \emph{Minkowski measurable (in dimension $D$)}.
More precisely, one defines the Minkowski dimension of the set $A$ as:
\begin{equation}
	\label{mink-dim}
	\dim_BA:=\inf\{r\geq 0|\mathcal{M}^r(A)=0\}=\sup\{r\geq 0|\mathcal{M}^r(A)=+\infty\}.
\end{equation}

In general, if the limit in \eqref{r-Mink} does not exist, one uses the notion of upper and lower $r$-dimensional Minkowski content by replacing the limit in \eqref{r-Mink} by the upper and lower limit, respectively, and denoting them by $\overline{\mathcal{M}}^r(A)$ and $\underline{\mathcal{M}}^r(A)$, respectively. 
Furthermore, one then defines the upper and lower Minkowski dimension analogously as in \eqref{mink-dim} and denotes them by $\overline{\dim}_BA$ and $\underline{\dim}_BA$, respectively.

Clearly, one always has that $\underline{\dim}_BA\leq\overline{\dim}_BA$ but note that even in the case when the Minkowski dimension exists, i.e., when $D=\underline{\dim}_BA=\overline{\dim}_BA$, the Minkowski content $\mathcal{M}^D(A)$, in general, needs not.
For our purpose, a weaker condition will be important, i.e., the case when $0<\underline{\mathcal{M}}^D(A)\leq\overline{\mathcal{M}}^D(A)<+\infty$.
One says that $A$ is \emph{Minkowski nondegenerate (in dimension $D$)} in this case.

A well known property of the Minkowski dimension is that it is invariant under bi-Lipschitz maps, i.e., if $f\colon A\to f(A)$ is bi-Lipschitz, then $\overline{\dim}_BA=\overline{\dim}_B(f(A))$ and analogously for the lower Minkowski dimension. 
Moreover, $A$ is Minkowski nondegenerate if and only if $f(A)$ is Minkowski nondegenerate.
On the other hand, note that Minkowski measurability is not preserved by bi-Lipschitz maps in general \cite{MeFr}, and the problem of characterizing Minkowski measurability is rather difficult and completely solved only for subsets of $\mathbb{R}$, see \cite{Falc95,LaPom,DKR,LaRaZu17,KomWin,LRZ-JFG}.
We refer the reader to, e.g., \cite{Falconer} for more details about the Minkowski dimension.
\smallskip

Let $(x_k)_{k\ge 0}$ and $(y_k)_{k\ge 0}$ be two sequences of positive real numbers. We say that $(x_k)_{k\ge 0}$ and $(y_k)_{k\ge 0}$ are comparable and write $x_k\simeq y_k$ as $k\to\infty$ if $1/\rho\le x_k/y_k\le \rho$ for some $\rho>0$ and for all $k\ge 0$. We use this notation in Section \ref{proof-theorem1} and Section \ref{section-applications}.

\section{Definitions and statement of results}
\label{section-statement}
In this section we state our main results.
We introduce the notions of the fractal sequence $\mathcal S$ of points on a smooth surface and associated chirp $\mathcal C_{\mathcal S}$, near a nilpotent contact point $p$ of even contact order and odd singularity order, and with finite slow divergence. We define the notion of the fractal codimension of $p$ in terms of the Minkowski dimension of $\mathcal S$ (see Definition \ref{def-intrin-codimen} in Section \ref{subsection-int-def-codim}). For the sake of readability, we first present this using a normal form for smooth
equivalence near $p$ (Section \ref{subsectionNintrinsic}). In the normal form coordinates, we give a complete fractal analysis of $p$ in terms of the Minkowski dimension of $\mathcal S$ (see Theorem \ref{theorem-1} in Section \ref{subsectionNintrinsic}). Theorem \ref{theorem-1} also gives the Minkowski dimension of $\mathcal C_{\mathcal S}$ if $p$ has finite (fractal) codimension. In Section \ref{subsectionintrinsic} we define the above notions (fractal sequence, chirp, fractal codimension) in an intrinsic way, using \cite{DDR-book-SF}, Section \ref{subsection-riemann} and Theorem \ref{theorem-1}. See Theorem \ref{theorem-Main} and Definition \ref{def-intrin-codimen} in Section \ref{subsection-int-def-codim}.

Assume that $p$ is a slow-fast Hopf point. We prove that in a classical Li\'{e}nard setting the codimension of $p$, defined in \cite{DRbirth}, and the fractal codimension of $p$ coincide (see Lemma \ref{lemma-trad-frac}). In Theorem \ref{theorem-Lienard} the fractal codimension of $p$ is connected with the cyclicity of $p$ (i.e. maximum number of limit cycles produced by $p$). We point out that we can intrinsically compute the fractal codimension of $p$, i.e. we don't need the classical Li\'{e}nard setting to determine it. 

\subsection{Calculating the Minkowski dimension in a normal form}
\label{subsectionNintrinsic}
Let $X_{\epsilon,\lambda}$ be a smooth family of vector fields on a smooth surface $M$ where $(\epsilon,\lambda)\sim (0,\lambda_0)\in \mathbb{R}\times\mathbb R^l$ and $\epsilon\ge 0$ is the  singular perturbation parameter. Assume that $X_{0,\lambda}$ has a curve $S_\lambda$ of singularities (we call $S_\lambda$ the critical curve) and that $p\in S_{\lambda_0}$ is a nilpotent singularity for $\lambda=\lambda_0$. Following \cite[Proposition 2.1]{DDR-book-SF}, there exists a local chart on $M$ around $p$ in which, after multiplication by
a smooth positive function, the system $X_{\epsilon,\lambda}$, for $(\epsilon,\lambda)\sim (0,\lambda_0)$, is given in
the following normal form:
\begin{equation}
\label{normal form 1}
    \begin{vf}
        \dot{x} &=& y-f(x,\lambda)  \\
        \dot{y} &=&\epsilon\left( g(x,\epsilon,\lambda)+\left(y-f(x,\lambda)\right)h(x,y,\epsilon,\lambda)\right),
    \end{vf}
\end{equation}
where  $f,g,h$ are smooth, $f(0,\lambda_0)=\frac{\partial f}{\partial x}(0,\lambda_0)=0$, $p=(0,0)$ in the local coordinates $(x,y)$ and where the overdot denotes differentiation with respect to the so-called fast time $t$. We have $S_\lambda=\{(x,y)|y=f(x,\lambda)\}$ in the coordinates $(x,y)$. When $\epsilon=0$, system \eqref{normal form 1} has horizontal fast movement away from $S_\lambda$. The contact order $c_{\lambda_0}(p)\ge 2$ is the order at $x=0$ of $f(x,\lambda_0)$ and the
singularity order $s_{\lambda_0}(p)\ge 0$ is the order at $x=0$ of $g(x,0,\lambda_0)$. Following \cite{DDR-book-SF}, the two notions of
order are independent of the chosen normal form. We suppose that $c_{\lambda_0}(p)$ and $s_{\lambda_0}(p)$ are finite and write 
\[f(x,\lambda_0)=x^{c_{\lambda_0}(p)}\tilde f(x)\] 
where $\tilde f$ is a smooth nowhere-zero function.
\textit{In the rest of this section we suppose that the parameter $\lambda$ in \eqref{normal form 1} is fixed: $\lambda=\lambda_0$.} 
\smallskip

\begin{figure}[htb]
	\begin{center}
		\includegraphics[width=6cm,height=2cm]{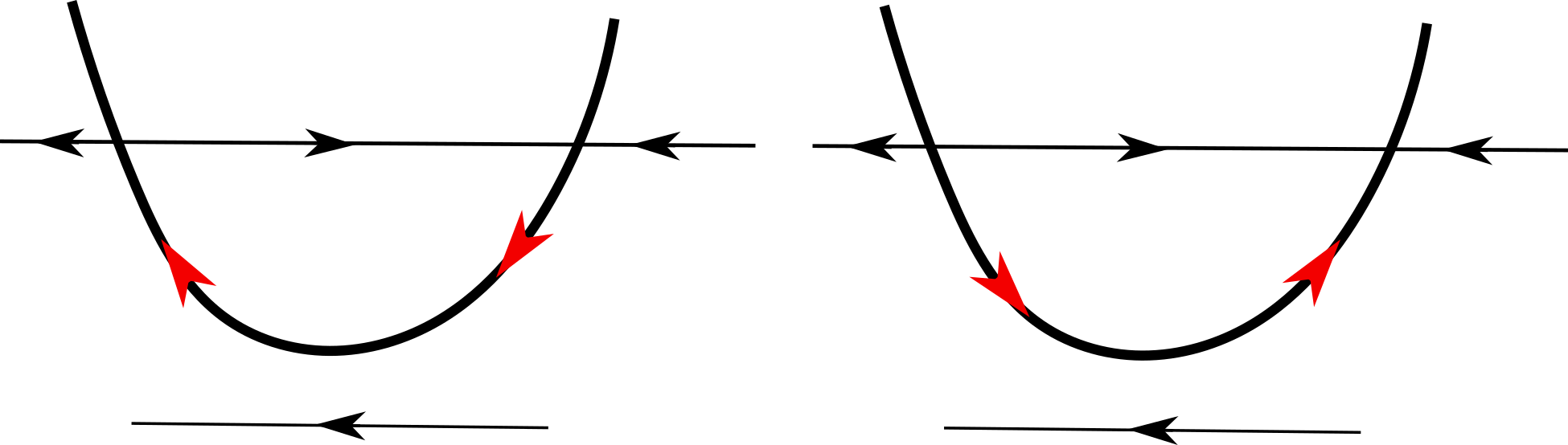}
		{\footnotesize
        \put(-139,-10){$(a)$}\put(-49,-10){$(b)$}
}
         \end{center}
	\caption{Fast dynamics of $\tilde X_\epsilon$ near the contact point with the indication of slow dynamics. (a) $g_m=-1$. (b) $g_m=+1$.}
	\label{fig-phase}
\end{figure}

 All the singularities on $S_{\lambda_0}$ are normally hyperbolic except for $p$ where we deal with a nilpotent singularity ($c_{\lambda_0}(p)\ge 2$ is finite!). Dynamics of \eqref{normal form 1} along the critical curve $S_{\lambda_0}$, for $\epsilon>0$, can be described using the following differential equation for so-called slow dynamics:
 \[\frac{dx}{ds}=\frac{g(x,0,\lambda_0)}{\frac{\partial f}{\partial x}(x,\lambda_0)}, \ \  x\ne 0,\]
 w.r.t. the so-called slow time $s=\epsilon t$, where $t$ is the fast time introduced in \eqref{normal form 1}.
 For more details see Section \ref{subsectionintrinsic} or \cite{DDR-book-SF}.
In \cite{DDR-book-SF} one can find a classification of all possible phase portraits of the limit $\epsilon=0$ in \eqref{normal form 1} with indication of the slow dynamics near $p$, depending on $c_{\lambda_0}(p)$ and $s_{\lambda_0}(p)$. In our slow-fast setting we have:\\
\\
\textbf{Assumption} We assume that:
\begin{itemize}
\item[\textbf{A}]{(Parity of $c_{\lambda_0}(p)$ and $s_{\lambda_0}(p)$)} $c_{\lambda_0}(p)$ is even and $s_{\lambda_0}(p)$ is odd. We write $n=c_{\lambda_0}(p)$ and $m=s_{\lambda_0}(p)$.
\item[\textbf{B}]{(Finite slow divergence)} $m\le 2(n-1)$.
\end{itemize}
From Assumption A it follows that: (a) $S_{\lambda_0}$ is a ``parabola-like'' critical curve where $p\in S_{\lambda_0}$ separates the normally attracting branch and the normally repelling branch of $S_{\lambda_0}$ and (b) the slow dynamics points from the attracting to the
repelling branch or from the repelling to the attracting branch (it cannot be directed toward $p$ or away from $p$ on both sides of $p$). If $\tilde{f}(0)>0$ (resp. $\tilde{f}(0)<0$), then the smooth diffeomorphism 
\[(x,y)\to (x \tilde f(x)^{\frac{1}{n}},y)\text{ }\Big(\text{resp. } (x,y)\to (-x  (-\tilde f(x))^{\frac{1}{n}},-y) \Big)\]
brings the system \eqref{normal form 1} ($\lambda=\lambda_0$), locally near $(x,y)=(0,0)$, into 
\begin{equation}
\label{normal form 2}
    \widetilde{X}_\epsilon: \ \begin{vf}
        \dot{x} &=& y-x^{n}  \\
        \dot{y} &=&\epsilon\left( g(x,\epsilon)+\left(y-x^n\right)h(x,y,\epsilon)\right),
    \end{vf}
\end{equation}
upon multiplication by a smooth strictly positive function, where $g,h$ are new smooth functions (we use the fact that $n$ is even). The two notions of order remain unchanged, i.e. the order at $x=0$ of $g(x,0)$ is equal to $m$. Using Assumption B we may assume that  the coefficient in front of $x^m$ in the expansion of $g(x,\epsilon)$ in powers of $x$ at $x=0$ is equal to $\pm 1$, after an $\epsilon$-depending rescaling in $(x,y,t)$. Thus, we can write the function $g(x,0)$ in \eqref{normal form 2} as 
    \[g(x,0)=g_m x^m+x^{m+1}\tilde{g}(x)\]
    where $g_m=\pm 1$ and $\tilde g$ is a smooth function. (As we will see later in this section, Assumption B will have another more important implication.) \textit{We have shown that, under Assumptions A and B, there exist smooth local coordinates $(x,y)$ in which (upon multiplication by a positive smooth function) the system $X_{\epsilon,\lambda_0}$ near $p\in M$, with $\epsilon\sim 0$, is given in the normal form $\widetilde{X}_\epsilon$.} In Fig. \ref{fig-phase} you can see the possible phase portraits of $\tilde{X}_0$, with indication of slow dynamics.
\smallskip

    The reason we want to use $\widetilde{X}_\epsilon$ is twofold. First, it is easy to define the notion of fractal codimension of $p$ in normal form $\widetilde{X}_\epsilon$ (see Definition \ref{def-codim} below). In Section \ref{subsectionintrinsic} (Definition \ref{def-intrin-codimen}) we define the notion of fractal codimension in a coordinate-free way. It will be clear that Definition \ref{def-codim} and Definition \ref{def-intrin-codimen} are equivalent in normal form coordinates \eqref{normal form 2}. Secondly, the normal form in \eqref{normal form 2} will be more convenient for computing the Minkowski dimension of fractal sequences and associated chirps introduced in Definition \ref{def-chirp} below.
    
\begin{definition}
\label{def-codim}
Consider $\widetilde{X}_\epsilon$ defined in \eqref{normal form 2}. We say that the contact point $p=(0,0)$ has finite (fractal) codimension $$\cdim(p):=j+1\ge 1$$ if 
\[\tilde{g}(x)+\tilde{g}(-x)=\alpha x^{2j}+O(x^{2j+2}), \ \alpha\ne 0.\]
If $j$ with the above property does not exist, we say that the codimension is infinite.
\end{definition}

In the normal form \eqref{normal form 2}, the critical curve $S_{\lambda_0}=\{y=x^n\}$ can be parameterized by the variable $x$. We define an open section $\sigma\subset\{x=0\}$ parameterized by the variable $y$ kept near $y=0$ ($y=0$ corresponds to the contact point $(0,0)\in\sigma$). For small and fixed $y,\tilde{y}>0$ ($y,\tilde y\in\sigma$) and for a sufficiently small $\rho>0$, we consider the slow divergence integral along $[-y^{1/n},-\rho]\cup [\rho,\tilde y^{1/n}]\subset S_{\lambda_0}$:
\begin{equation}\label{SDI-between}I(y,\tilde y,\rho)=-\int_{-y^{1/n}}^{-\rho}\frac{1}{g(x,0)}(nx^{n-1})^2dx-\int_{\rho}^{\tilde y^{1/n}}\frac{1}{g(x,0)}(nx^{n-1})^2dx.
\end{equation}
The slow divergence of $\tilde X_0$ is given by $-nx^{n-1}$ and $ds=\frac{nx^{n-1}}{g(x,0)}dx$ where $s$ is the slow time. For more details see \cite{DDR-book-SF} or Section \ref{subsectionintrinsic}.
Assumption B implies that the following limit, which represents the slow divergence integral along $[-y^{1/n},\tilde y^{1/n}]\subset S_{\lambda_0}$ is finite:
\begin{equation}
    \label{SDI-contact}
    I(y,\tilde y)=\lim_{\rho\to 0^+}I(y,\tilde y,\rho)=-\int_{-y^{1/n}}^{\tilde y^{1/n}}\frac{1}{g(x,0)}(nx^{n-1})^2dx.
\end{equation}
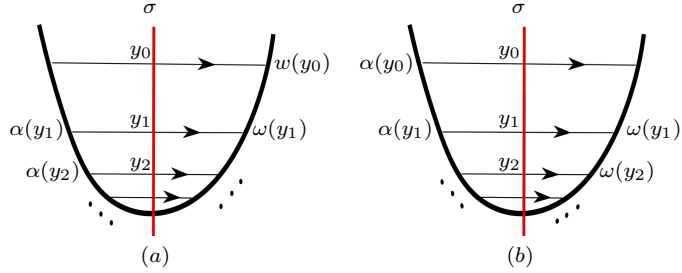
\begin{figure}[htb]
	\begin{center}
		\includegraphics[width=8cm,height=2.8cm]{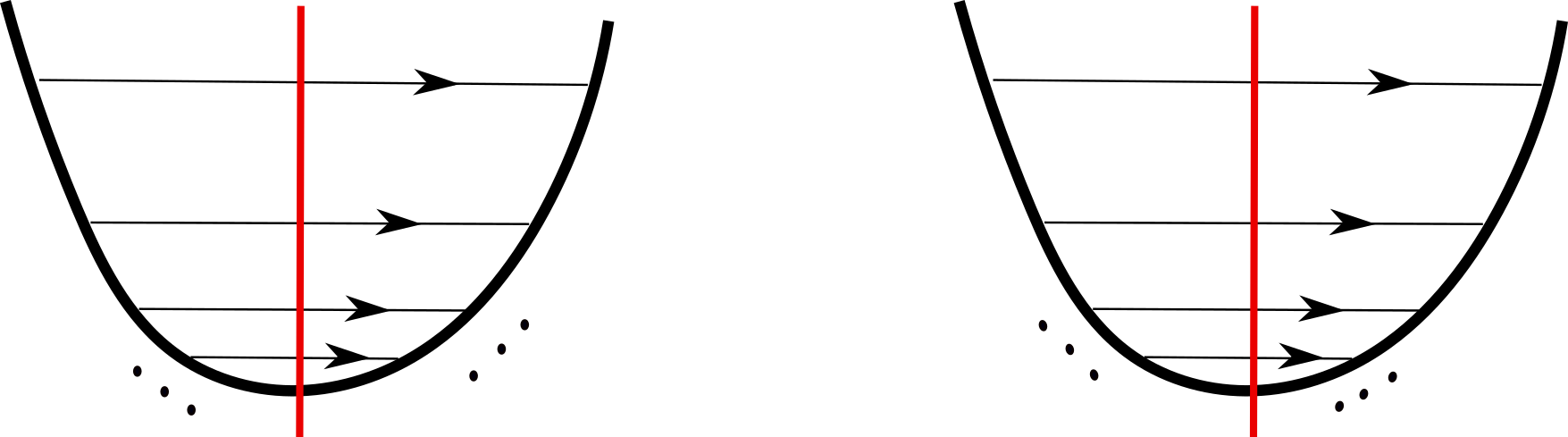}
		{\footnotesize
        \put(-139,63){$w(y_0)$}\put(-107,63){$\alpha(y_0)$}
        \put(-238,38){$\alpha(y_1)$}
        \put(-147,38){$\omega(y_1)$}
        \put(-232,22){$\alpha(y_2)$}
        \put(-100,38){$\alpha(y_1)$}
        \put(-7,38){$\omega(y_1)$}
        \put(-17,22){$\omega(y_2)$}
        \put(-193,68){$y_0$}
        \put(-193,42){$y_1$}
        \put(-193,26){$y_2$}
        \put(-188,84){$\sigma$}
        \put(-55,68){$y_0$}
        \put(-55,42){$y_1$}
        \put(-55,26){$y_2$}
        \put(-50,84){$\sigma$}
        \put(-189,-10){$(a)$}\put(-52,-10){$(b)$}
        }

         \end{center}
	\caption{A fractal sequence starting at $(0,y_0)$ defined near the contact point $(x,y)=(0,0)$ where $\alpha(y)=\{(-y^{1/n},y)\}$ is the $\alpha$-limit of the fast orbit through $y\in \sigma$ and $\omega(y)=\{(y^{1/n},y)\}$ is the $\omega$-limit of the same orbit. (a) We use $I(y_{k+1},y_{k})=0$ to generate $(y_k)_{k\ge 0}$. (b) We use $I(y_{k},y_{k+1})=0$ to generate $(y_k)_{k\ge 0}$.}
	\label{fig-sequence}
\end{figure}
Since $2(n-1)-m$ is odd (Assumption A), the integrand in \eqref{SDI-contact} changes sign as $x$ passes through $x=0$, and for each $y>0$ (resp. $\tilde y>0$) small enough we can therefore find a unique $\tilde y>0$ (resp. $y>0$) such that $I(y,\tilde y)=0$. In order to avoid the possibility that $I(y,y)=0$ for some $y>0$ and $y\sim 0$, we make the following assumption concerning the divergence integral $I(y):=I(y,y)$ (clearly, $I(0)=0$):\\
\\
\textbf{Assumption C} We assume that $y=0$ is an isolated zero of $I(y)$, meaning that there is a small $y^*>0$ such that $I$ is nonzero on $]0,y^*[$.\\
\\
Since $I$ is continuous, from Assumption C it follows that $I$ is either positive or negative on $]0,y^*[$. Let $y_0\in ]0,y^*[$ be arbitrary and fixed. We distinguish two possibilities to generate a strictly decreasing sequence $(y_k)_{k\ge 0}$, with the initial point $y_0$, converging to zero using the ``entry-exit'' relation $I(y,\tilde y)=0$:
\begin{itemize}
\item ($I(y)>0$ on $]0,y^*[$ and $g_m=+1$) or ($I(y)<0$ on $]0,y^*[$ and $g_m=-1$) Here the sequence $(y_k)_{k\ge 0}$, defined by $I(y_{k+1},y_{k})=0$ with $k\ge 0$, is decreasing and converges to zero. See Fig. \ref{fig-sequence}(a).
    \item ($I(y)>0$ on $]0,y^*[$ and $g_m=-1$) or ($I(y)<0$ on $]0,y^*[$ and $g_m=+1$) In this case the sequence $(y_k)_{k\ge 0}$, defined by $I(y_k,y_{k+1})=0$ with $k\ge 0$, is decreasing and converges to zero. See Fig. \ref{fig-sequence}(b).
\end{itemize}
\begin{remark}
\label{remark-alpha}
In Section \ref{proof-theorem1} we will prove that, if $\tilde X_{\epsilon}$ has finite codimension $\cdim(p)=j+1\ge 1$, then $I$ from Assumption C is positive (resp. negative) if and only if $\alpha>0$ (resp. $\alpha<0$) where $\alpha$ is introduced in Definition \ref{def-codim}.
\end{remark}
Now we introduce the notion of the fractal sequence and the associated chirp.
\begin{definition} Consider $\widetilde{X}_\epsilon$ defined in \eqref{normal form 2}. Let $y_0\in ]0,y^*[$ be arbitrary and fixed and let $(y_k)_{k\ge 0}$ be the (unique) sequence defined above. We call the set $\mathcal S=\{(0,y_k) \ | \ k\ge 0\}$ the fractal sequence of the contact point \eqref{normal form 2} with the initial point $(0,y_0)$, and the set 
\[\mathcal{C}_{\mathcal S}=\cup_{k=0}^{\infty} \mathcal{C}_k\subset \mathbb R^2, \]
with $\mathcal{C}_k=]-y_k^{1/n},y_k^{1/n}[\times\{y_k\}$, the associated chirp.
\label{def-chirp}
\end{definition}
Finally we state the main result of this section.
\begin{theorem}
\label{theorem-1}
Consider a smooth slow-fast system $\tilde X_{\epsilon}$ defined in \eqref{normal form 2}, that satisfies Assumptions A--C. Suppose that $\tilde X_{\epsilon}$ has finite fractal codimension $\cdim(p)=j+1\ge 1$. Then $\mathcal S$ and $\mathcal{C}_{\mathcal S}$, given in Definition \ref{def-chirp}, are Minkowski nondegenerate and  \[\dim_B\mathcal S=\frac{2j+1}{n+2j+1}\in ]0,1[\text{  and  }\dim_B\mathcal{C}_{\mathcal S}=\frac{n+4j+1}{n+2j+1}\in [1,2[.\]
Moreover, when the codimension is infinite, we have $\dim_B\mathcal S=1$. The results do not depend on the choice of the initial point $y_0\in ]0,y^*[$.
\end{theorem}
Theorem \ref{theorem-1} will be proved in Section \ref{proof-theorem1}.
\begin{remark}\label{remark-disc-val}
1. Theorem \ref{theorem-1} implies that, for a fixed contact order $n$ and singularity order $m$ in $\tilde X_{\epsilon}$, there exists a 1-1 correspondence between $\dim_B\mathcal S$ and the fractal codimension $\cdim(p)=j+1\in\mathbb{N}\cup \{+\infty\}$ of $\tilde X_{\epsilon}$. Since we can intrinsically define the notion of contact order and $\dim_B\mathcal S$ (see Section \ref{subsectionintrinsic}), this observation helps us to intrinsically define the notion of fractal codimension (see Definition \ref{def-intrin-codimen}). 

2. Following Theorem \ref{theorem-1}, if $n=2$ (and hence $m=1$), then $\dim_B\mathcal S$ can take the following discrete set of values: $\frac{1}{3},\frac{3}{5},\frac{5}{7},\dots,1$.

3. There also exists a bijective correspondence between $\dim_B\mathcal{C}_{\mathcal S}$ and finite codimension $\cdim(p)=j+1\in\mathbb{N}$ of $\tilde X_{\epsilon}$.
We believe that $\dim_B\mathcal{C}_{\mathcal S}=2$ when $\tilde X_{\epsilon}$ has infinite fractal codimension. This is a topic of further study.

4. Note that $\dim_B\mathcal S$ and $\dim_B\mathcal{C}_{\mathcal S}$ given in Theorem \ref{theorem-1} are independent of the singularity order $m$.
\end{remark}

\subsection{Intrinsic nature of fractal codimension of contact points}
\label{subsectionintrinsic}
In Section \ref{subsub-general} we recall some basic definitions from \cite{DDR-book-SF} for smooth slow-fast systems defined on a smooth surface. We intrinsically define the notions of critical curve, fast foliation, contact order, singularity order, slow vector field, slow divergence integral, etc. In Section \ref{subsection-int-def-codim} we (intrinsically) define the notion of fractal codimension and state our main results.

\subsubsection{Slow-fast model and basic definitions}\label{subsub-general}
  Let $X_{\epsilon,\lambda}$ be a smooth $\epsilon-$slow-fast system on a smooth surface $M$. We write \begin{equation}\label{system-generalni} X_{\epsilon,\lambda}=X_{0,\lambda}+\epsilon Q_\lambda+O(\epsilon^2)
  \end{equation}
where $X_{0,\lambda}$ has a set of non-isolated singularities $S_\lambda$ for each $\lambda\sim\lambda_0$. We suppose that for every singularity $p\in S_{\lambda_0}$  there is an open neighborhood $V\subset M$ of $p$ such that $X_{0,\lambda}=F_\lambda Z_\lambda$ on $V$, where $F_\lambda$ (resp. $Z_\lambda$) is a smooth family of  functions (resp. vector fields without singularities for each $\lambda\sim\lambda_0$). We further assume that $dF_\lambda(p)\ne 0$, $\forall p\in\{F_\lambda=0\}$. We call $\{V,Z_\lambda,F_\lambda\}$ an admissible expression for $X_{0,\lambda}$ near $p$ (see \cite{DDR-book-SF}). For example, in \eqref{normal form 1} (system \eqref{normal form 2} is a special case of \eqref{normal form 1}) we have 
\[Z_\lambda=\frac{\partial}{\partial x}, \ F_\lambda(x,y)=y-f(x,\lambda), \ Q_\lambda=\left(g(x,0,\lambda)+F_\lambda(x,y)h(x,y,0,\lambda)\right)\frac{\partial}{\partial y}. \]
Note that the pair $(Z_\lambda,F_\lambda)$ is not unique (we can replace $Z_\lambda$ with $\xi_\lambda Z_\lambda$ and $F_\lambda$ with $\frac{1}{\xi_\lambda} F_\lambda$ where $\xi_\lambda$ is a nowhere zero function). We denote by $t$ so-called fast time related to \eqref{system-generalni} (for more details see \cite{DDR-book-SF}).
\smallskip

 Clearly, $V\cap S_\lambda=\{F_\lambda=0\}$ and $S_\lambda$ is a $1$-dimensional submanifold of $M$ (often called the critical curve). We denote by $\mathcal F_\lambda$ a smooth $1$-dimensional foliation on $M$ tangent to $Z_\lambda$ in each admissible local expression $\{V,Z_\lambda,F_\lambda\}$ for $X_{0,\lambda}$. We call $\mathcal F_\lambda$ the fast foliation of $X_{0,\lambda}$, and the orbits of the fast flow of $X_{0,\lambda}$, away from $S_\lambda$, are located inside the leaves of the fast foliation (the leaf
through $p\in M$ is denoted by $l_{\lambda,p}$). It is clear that the notions of critical curve, fast foliation, fast flow, etc. are defined in an intrinsic (i.e. a coordinate-free) way. For more details see \cite{DDR-book-SF}. In \eqref{normal form 1}, the fast foliation is given by horizontal lines.
\smallskip

If the linear part of $X_{0,\lambda}$ at $p\in S_\lambda$ has an eigenvalue $E_\lambda(p)$ different from zero, then we say that $p$ is normally hyperbolic (attracting when $E_\lambda(p)<0$ or repelling when $E_\lambda(p)>0$). (There is always one zero eigenvalue with eigenspace $T_pS_\lambda$.) The nonzero eigenvalue $E_\lambda(p)$ is also intrinsically defined: its eigenspace is $T_pl_{\lambda,p}$ and $E_\lambda(p)$ is equal to the trace of $DX_{0,\lambda}(p)$. When the linear part of $X_{0,\lambda}$ at $p$ has both eigenvalues equal to zero, the point $p$ is called a contact point (between $S_\lambda$ and $\mathcal F_\lambda$). The above mentioned conditions on $F_\lambda$ and $Z_\lambda$ imply that all contact points are nilpotent. In this paper we deal with an isolated (nilpotent) contact point. 
\smallskip

Following \cite{DDR-book-SF}, we can use the following intrinsic definition of contact order and singularity order of a contact point $p$ at level $\lambda=\lambda_0$: the contact order $n=n_{\lambda_0}(p)$ is the contact at $p$ between $S_{\lambda_0}$ and the leaf $l_{\lambda_0,p}$, and, for any admissible expression  $\{V,Z_\lambda,F_\lambda\}$ for $X_{0,\lambda}$ near $p$ and for any area form $\Omega$ on $V$, the singularity order $m=s_{\lambda_0}(p)$ of $p$ is the order at $p$ of the function $\Omega(Q_{\lambda_0},Z_{\lambda_0})|_{S_{\lambda_0}\cap V}: \tilde p\in S_{\lambda_0}\cap V\mapsto \Omega(Q_{\lambda_0},Z_{\lambda_0})(\tilde p)$. From \cite[Lemma 2.1]{DDR-book-SF} it follows that the definition does not depend on the choice of the admissible expression near $p$ and $\Omega$. This intrinsic definition coincides with the definition given in Section \ref{subsectionNintrinsic}, using normal form coordinates \eqref{normal form 1} for smooth equivalence (see \cite{DDR-book-SF}). \textit{We suppose that the contact point $p$ for $\lambda=\lambda_0$ satisfies Assumptions A and B defined in Section \ref{subsectionNintrinsic}}.
\smallskip

The notion of slow vector field on normally hyperbolic portions of $S_\lambda$ is also defined in an intrinsic way. If $p\in S_\lambda$ is normally hyperbolic, then we denote by $Q^*_\lambda(p)\in T_pS_\lambda$ the linear projection of $Q_\lambda(p)$ on $T_pS_\lambda$ in the direction of the transverse eigenspace $T_pl_{\lambda,p}$. We call $Q^*_\lambda$ the slow vector field and its flow the slow dynamics. The time variable of the slow dynamics is the 
slow time $s=\epsilon t$ where $t$ is the fast time attached to $X_{\epsilon,\lambda}$. This definition is equivalent to the classical one using center manifolds (see e.g. \cite{DDR-book-SF}). In \eqref{normal form 1} we have 
\[Q^*_\lambda=\frac{g(x,0,\lambda)}{\frac{\partial f}{\partial x}(x,\lambda)}\frac{\partial}{\partial x}+g(x,0,\lambda)\frac{\partial}{\partial y}\]
when $\frac{\partial f}{\partial x}(x,\lambda)\ne 0$. 
\smallskip

We give now an intrinsic definition of the notion of slow divergence integral. If $\gamma_\lambda\subset S_\lambda$ is a segment containing normally hyperbolic singularities different from zeros of the slow vector field $Q^*_\lambda$, then the slow divergence integral along $\gamma_\lambda$ is defined by
\begin{equation}
\label{SDI-intr}
I(\gamma_\lambda):=\int_{\gamma_\lambda}E_\lambda ds.
\end{equation}
We can compute $I(\gamma_\lambda)$ using a parameterization of $\gamma_\lambda$. More precisely, if $p_\lambda:[r_1,r_2]\to M$ is a parameterization of the segment $\gamma_\lambda$ and $Q_\lambda^*(p_\lambda(r))=q_\lambda(r)\frac{dp_{\lambda}}{dr}$, then \eqref{SDI-intr} becomes
\begin{equation}
\label{SDI-intr1}
I(\gamma_\lambda)=I(r_1,r_2)=\int_{r_1}^{r_2}E_\lambda(p_\lambda(r))\frac{dr}{q_\lambda(r)}\nonumber
\end{equation}
where $q_\lambda$ has no zeros on $[r_1,r_2]$ because $Q^*_\lambda$ has no zeros on $\gamma_\lambda$ (see \cite[Proposition 5.3]{DDR-book-SF}). For example,
in \eqref{normal form 1} we can take (see the expression for $Q^*_\lambda$):
\[r=x, \ p_\lambda (x)=(x,f(x,\lambda)), \ q_\lambda(x)=\frac{g(x,0,\lambda)}{\frac{\partial f}{\partial x}(x,\lambda)}, \ E_\lambda(p_\lambda(x))=-\frac{\partial f}{\partial x}(x,\lambda),\]
and, therefore, we obtain
\begin{equation}\label{SDI-intr2}I(x_1,x_2)=-\int_{x_1}^{x_2}\frac{\left(\frac{\partial f}{\partial x}(x,\lambda)\right)^2}{g(x,0,\lambda)}dx.
\end{equation}
Note that the expression in \eqref{SDI-between} comes from \eqref{SDI-intr2} where $f(x,\lambda)=x^n$ and we replace $g(x,0,\lambda)$ with $g(x,0)$. 
\smallskip

\textit{One important property of the slow divergence integral \eqref{SDI-intr} we will use throughout this paper is its invariance under changes of coordinates, time
reparameterizations and changes of the area form on $M$.} 

\subsubsection{Fractal codimension and statement of results}
\label{subsection-int-def-codim}
Let's fix $\lambda=\lambda_0$.
Suppose that $p\in S_{\lambda_0}$ is a (nilpotent) contact point of $X_{0,\lambda_0}$ and take a small neighborhood $V$ of $p$. From Assumptions A and B  it follows that we can compute  the slow divergence integral \eqref{SDI-intr} along any segment $\gamma\subset S_{\lambda_0}$, $\gamma\subset V$, containing $p$ in its interior. (The segment $\gamma$ consists of a normally attracting branch, a normally repelling branch and $p$ between them.) Indeed, since the notion of slow divergence integral does not depend on the local coordinate system and does not change under time reparameterizations, it suffices to compute the integral along $\gamma$ in a normal form \eqref{normal form 2} for smooth equivalence. The expression of the integral along $\gamma$ in the normal form  \eqref{normal form 2} is given in \eqref{SDI-contact} and it is finite under Assumption B (see Section \ref{subsectionNintrinsic}).
\smallskip

\begin{figure}[htb]
	\begin{center}
		\includegraphics[width=6cm,height=3.3cm]{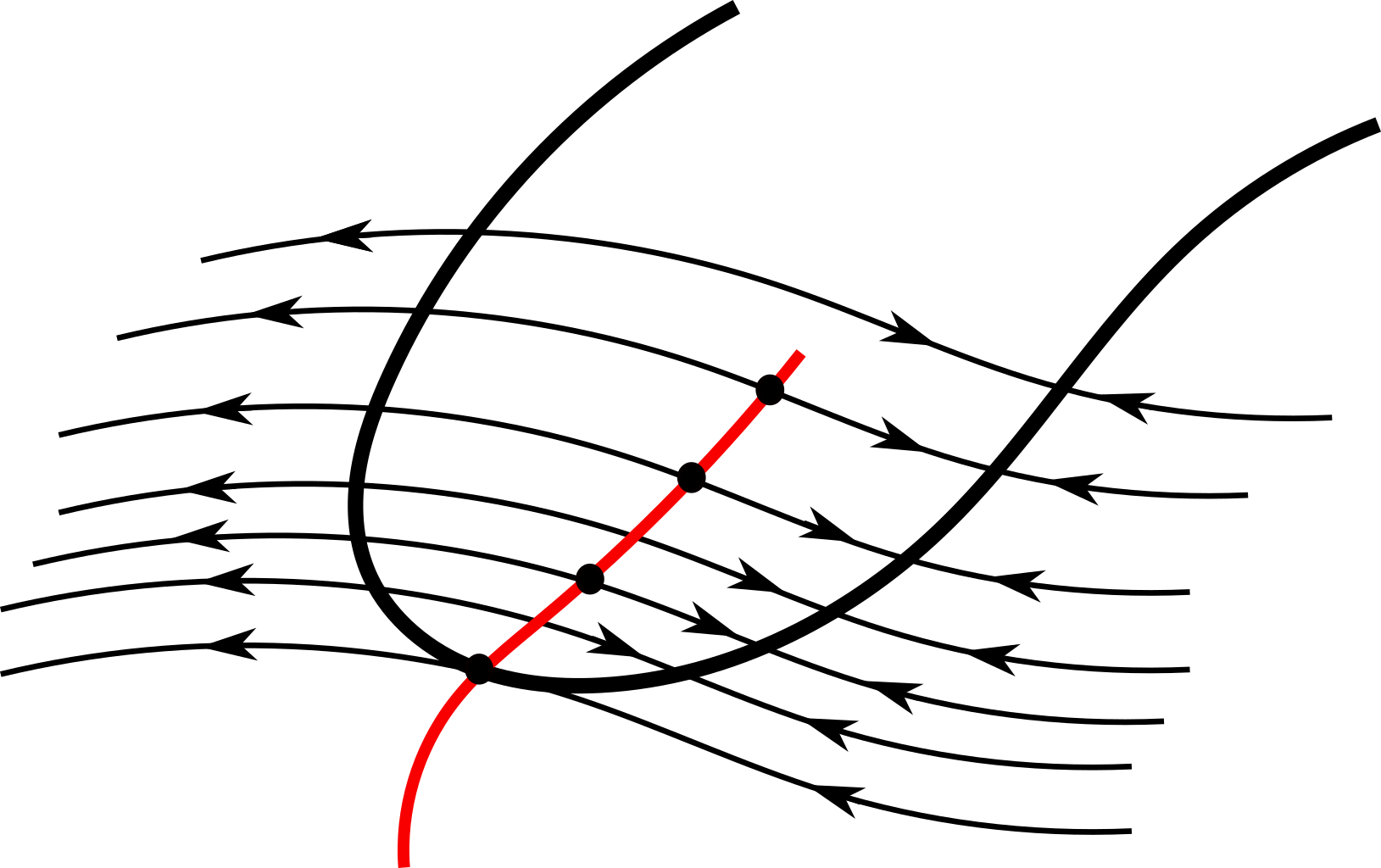}
		{\footnotesize
        \put(-122,-7){$\sigma$}
        \put(-113,15){$p$}
\put(-83,56){$p_0$}\put(-91,46){$p_1$}\put(-107,34){$p_2$}\put(-130,8){$\sigma_2$} \put(-80,43){$\sigma_1$} \put(-20,82){$S_{\lambda_0}$}\put(-143,52){$\alpha(p_1)$}\put(-54,33){$\omega(p_1)$}
}
         \end{center}
	\caption{A section $\sigma=\sigma_1\cup\{p\}\cup \sigma_2$ transverse to the fast foliation $\mathcal F_{\lambda_0}$, with $\sigma_1$ located ``inside'' the critical curve $S_{\lambda_0}$ and a fractal sequence $(p_k)_{k\ge 0}$ monotonically tending to $p$.}
	\label{fig-onM}
\end{figure}
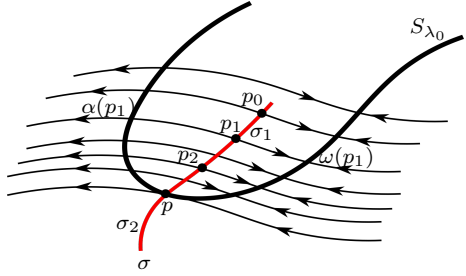

Let $\sigma\subset M$ be a smooth regular (open) section transverse to
 the fast foliation $\mathcal{F}_{\lambda_0}$ such that $p\in \sigma$. We write $\sigma=\sigma_1\cup\{p\}\cup \sigma_2$ where $\sigma_1$ is located in the region bounded by the ``parabola-like'' critical curve $S_{\lambda_0}$ (see Fig. \ref{fig-onM}). For $\tilde p\in\sigma_1$ we denote by $\alpha(\tilde p)\in S_{\lambda_0}$ (resp. $\omega(\tilde p)\in S_{\lambda_0}$) the $\alpha$-limit point (resp. the $\omega$-limit point) of the fast orbit of $X_{0,\lambda_0}$ through $\tilde p$. The singularity $\alpha(\tilde p)$ (resp. $\omega(\tilde p)$) is normally repelling (resp. attracting). Further, we denote by $I(\tilde p)$, with $\tilde p\in\sigma_1$, the slow divergence integral along the segment $\gamma\subset S_{\lambda_0}$ from $\alpha(\tilde p)$ to $\omega(\tilde p)$, and by $I(\tilde p,\bar p)$, with $\tilde p,\bar p\in\sigma_1$, the slow divergence integral along $\gamma\subset S_{\lambda_0}$ from $\alpha(\tilde p)$ to $\omega(\bar p)$. Clearly, we have $I(\tilde p)=I(\tilde p,\tilde p)$ and $I(\tilde p)\to 0$ as $\tilde p\to p$.
Now, using the invariance of the slow divergence integral under smooth equivalences, Assumption C from Section \ref{subsectionNintrinsic} becomes:\\
\\
\textbf{Assumption C'} We assume that $I(\tilde p)\ne 0$ for all $\tilde p\in\sigma_1$. 
\smallskip

\begin{remark}
Note that Assumption C' does not depend on the choice of section $\sigma$. If we take another section $\sigma'=\sigma_1'\cup\{p\}\cup \sigma_2'$ in a sufficiently small neighborhood of $p$, with the same property as $\sigma$, then $I(\tilde p)\ne 0$ for all $\tilde p\in\sigma_1'$. This follows from the definition of $I(\tilde p)$ and the fact that each fast orbit of $X_{0,\lambda_0}$ that intersects $\sigma_1'$ also intersects $\sigma_1$.
\end{remark}

Let $p_0\in\sigma_1$ be arbitrary and fixed. From Section \ref{subsectionNintrinsic} and Assumption C' we know that precisely one of the recursive formulas, $I(p_k,p_{k+1})=0$, $k\ge 0$, or $I(p_{k+1},p_{k})=0$, $k\ge 0$, generates  a sequence $(p_k)_{k\ge 0}$ in $\sigma_1$ that monotonically converges to the contact point $p$. We can identify the formula which generates such a sequence in the following coordinate free way:
\begin{itemize}
    \item ($I(\tilde p)>0$, for $\tilde p\in\sigma_1$, and the slow vector field $Q^*_{\lambda_0}$ points from the repelling part to the attracting part of $S_{\lambda_0}$) or ($I(\tilde p)<0$, for $\tilde p\in\sigma_1$, and $Q^*_{\lambda_0}$ points from the attracting part to the repelling part of $S_{\lambda_0}$) In this case for any fixed $p_0\in\sigma_1$ the (unique) sequence $(p_k)_{k\ge 0}$, generated by $I(p_{k+1},p_{k})=0$, $k\ge 0$, monotonically converges to $p$. 
    
    \item ($I(\tilde p)>0$, for $\tilde p\in\sigma_1$, and $Q^*_{\lambda_0}$ points from the attracting part to the repelling part of $S_{\lambda_0}$) or ($I(\tilde p)<0$, for $\tilde p\in\sigma_1$, and $Q^*_{\lambda_0}$ points from the repelling part to the attracting part of $S_{\lambda_0}$) In this case we use $I(p_k,p_{k+1})=0$, $k\ge 0$, in order to produce such a sequence, for any $p_0\in\sigma_1$. 
\end{itemize}

The set $\mathcal S=\{p_k \ | \ k\ge 0\}$ is called the fractal sequence starting at $p_0$, and the set $\mathcal C_{\mathcal S}=\cup_{k=0}^\infty \mathcal C_k$, where $\mathcal C_k$ is the fast orbit of $X_{0,\lambda_0}$ through $p_k$, is the associated chirp. It is clear that these two notions are defined in a coordinate-free way.
\begin{remark}
\label{remarkSS'}
Let $\mathcal S=\{p_k \ | \ k\ge 0\}$ be a fractal sequence on $\sigma$. If we choose another section $\sigma'=\sigma_1'\cup\{p\}\cup \sigma_2'$ with the same property as $\sigma$, then we can define a set $\mathcal S '=\{p_k' \ | \ k\ge 0\}\subset \sigma'_1$ where $p_k'$ is the (unique) point at which the fast orbit of $X_{0,\lambda_0}$ through $p_k$ hits $\sigma_1'$ (we assume that the fast orbit through $p_0$ intersects $\sigma_1'$). From the definition of $I(\tilde p,\bar p)$ and the fact that $p_k$ and $p_k'$ have the same $\alpha$- and $\omega$-limit points, it follows that $\mathcal S'$ is a fractal sequence, starting at $p_0'$, generated by the same recursive formula as $\mathcal S$. It is clear that $\mathcal C_{\mathcal S}=\mathcal C_{\mathcal S'}$.
\end{remark}

\begin{theorem}
\label{theorem-Main}
Consider a smooth slow-fast system $X_{\epsilon,\lambda}$ on $M$, given in \eqref{system-generalni}. Let $p\in S_{\lambda_0}$ be a (nilpotent) contact point of $X_{0,\lambda_0}$ of contact order $n$ and singularity order $m$ that satisfies Assumptions A, B and C', and let $\mathcal S$ be a fractal sequence defined above. Then $\dim_B\mathcal S$ exists and \begin{equation}\label{equations-bijection}
  \dim_B\mathcal S\in \left\{\frac{2j+1}{n+2j+1} \ | \ j\in\mathbb N_0 \right\}\cup \{1\}.  
\end{equation} 
If $\dim_B\mathcal S\ne 1$, then $\mathcal S $ and $\mathcal C_{\mathcal S}$ are Minkowski nondegenerate. Furthermore, the Minkowski dimension and Minkowski nondegeneracy of $\mathcal S$ are coordinate free notions which do not depend on the choice of the section $\sigma$, transversally
cutting the leaf $l_{\lambda_0,p}$ at the contact point $p$, the first element $p_0$ of the sequence $(p_k)_{k\ge 0}$ from $\mathcal S$, and the metric on $M$. The same is true for chirps if $\dim_B\mathcal S\ne 1$.
\end{theorem}
We prove Theorem \ref{theorem-Main} in Section \ref{section-proof-Main}. The reason why Theorem \ref{theorem-Main} is important is twofold. First, we are allowed to compute $\dim_B\mathcal S $ in any local coordinates $(x,y)$, e.g. in the normal form coordinates \eqref{normal form 2} using the Euclidean metric. On the other hand, since for fixed $n$ and $m$ there exists a bijective correspondence between $j\in \mathbb{N}_0\cup \{+\infty\}$ and $\dim_B\mathcal S$ (see Remark \ref{remark-disc-val}), we have the following natural definition of the (fractal) codimension of the contact point $p$.
\begin{definition}\label{def-intrin-codimen}
Let $X_{\epsilon,\lambda}$ be a slow-fast system as in Theorem \ref{theorem-Main} (i.e. satisfying the same conditions). If $\dim_B\mathcal S<1$, we say that the contact point $p$ for $\lambda=\lambda_0$ has finite fractal codimension $$\cdim(p):=j+1\ge 1$$
where 
\[j=\frac{(n+1)\dim_B\mathcal S-1}{2(1-\dim_B\mathcal S)}\in\mathbb{N}_0.\]
If $\dim_B\mathcal S=1$, then we say that the fractal codimension of $p$ is infinite.
\end{definition}
Note that Theorem \ref{theorem-1} implies that the above intrinsic definition of the codimension coincides with the definition of the codimension in normal form coordinates \eqref{normal form 2} (see Definition \ref{def-codim}).
\smallskip

In the rest of this section we focus on a contact point $p$ at $\lambda=\lambda_0$ of Morse type (i.e. of contact order $2$) and singularity order $1$. It is clear that Assumptions A and B are satisfied in this case. Since the singularity order of $p$ is $1$, we can define the notion of singularity index (see \cite{DDR-book-SF}): $p$ has singularity index $\pm 1$ when $\text{sign} \frac{\partial g}{\partial x}(0,0,\lambda_0)=\mp 1$ where $g$ comes from any normal form \eqref{normal form 1}. (For an intrinsic definition of the singularity index we refer to \cite[Lemma 2.3]{DDR-book-SF}.) For $\epsilon>0$ and $\lambda\sim\lambda_0$, $X_{\epsilon,\lambda}$ has a hyperbolic saddle near $p$ if $p$ has singularity index $-1$ or a center/focus/node if $p$ has singularity index $+1$.

\smallskip
When $p$ has singularity index $+1$, then we say that $p$ is a slow-fast Hopf point. Such a slow-fast Hopf point can produce limit cycles after perturbation.  We say that the cyclicity of a slow fast Hopf point $p$ in $X_{\epsilon,\lambda}$ is bounded by $N\in\mathbb{N}_0$ if we can find $\epsilon_0>0$, a neighborhood $U$ of $\lambda_0$ in the $\lambda$-space and a neighborhood $V$ of $p$ such that $X_{\epsilon,\lambda}$ has at most $N$ limit cycles in $V$ for each $(\epsilon,\lambda)\in [0,\epsilon_0]\times U$. We call the smallest $N$ with this property the cyclicity of $p$ in the slow-fast family $X_{\epsilon,\lambda}$ and denote it by $\cycl(X_{\epsilon,\lambda},p)$.

\smallskip
Consider a Li\'{e}nard system with a slow-fast Hopf point at $(x,y)=(0,0)$:
\begin{equation}
\label{model-Hopf-for-def}
    \begin{vf}
        \dot{x} &=& y-x^2+x^3h_1(x,\epsilon,\lambda)  \\
        \dot{y} &=&\epsilon\big(b(\epsilon,\lambda)-x\big),
    \end{vf}
\end{equation}
where $h_1,b$ are smooth functions and $b(0,\lambda_0)=0$. Following \cite[Definition 1.4]{DRbirth}, the system \eqref{model-Hopf-for-def} has a finite codimension at $(\epsilon,\lambda)=(0,\lambda_0)$ if there exists $j\in\mathbb N_0$ such that  
\begin{equation}\label{h1even}
h_1(x,0,\lambda_0)+h_1(-x,0,\lambda_0)=\alpha x^{2j}+O(x^{2j+2}), \ \ \alpha\ne 0.
\end{equation}
If such a $j\in\mathbb N_{0}$ exists, then we say that the slow-fast Hopf point has a codimension equal to $j+1\ge 1$. If $j$ with the property \eqref{h1even} does not exist, then we say that the codimension is infinite.

\begin{remark}
\label{remark-conjecture}
If the smooth Li\'{e}nard system \eqref{model-Hopf-for-def} has finite codimension $j+1$, then the cyclicity of $(x,y)=(0,0)$ in \eqref{model-Hopf-for-def} is finite. Furthermore, if the conjecture stated below holds for all $i\le j+1$, then the cyclicity of $(x,y)=(0,0)$ in \eqref{model-Hopf-for-def} is bounded by $j+1$ (see \cite[Theorem 7.10]{DDR-book-SF} or \cite{DRbirth}). Let $H(x,y)=e^{-2y}\left(y-x^2+\frac{1}{2}\right)$ and let $\gamma_h$ be the oval described by the level set $\{H(x,y)=h\}$, where $h\in ]0,\frac{1}{2}[$. When $h=\frac{1}{2}$, the level set is the origin. We consider the integrals $J_i(h)=\int_{\gamma_h}x^{2i-1}dy$, $h\in]0,\frac{1}{2}]$ and $i=0,1,2,\dots$. The conjecture from \cite{DRbirth} can now be stated as follows: for each $i\ge 0$, the ordered set $(J_0,\dots, J_i) $ forms a strict Chebyshev system on $[h_0,\frac{1}{2}]$ for any $h_0\in ]0,\frac{1}{2}[$. Strict Chebyshev systems are defined in \cite{DRbirth}. The conjecture was proven in \cite{Villa} or \cite{Li-Liu-Llibre} for $i\le 2$. The
conjecture, for any $i\ge 0$, has been solved recently by Chengzhi Li and Changjian Liu in \cite{Li-conjecture-solved}.

If \eqref{model-Hopf-for-def} is an analytic slow-fast family, then the cyclicity of $(x,y)=(0,0)$ is finite (see \cite{DRbirth}).

\end{remark}

The above ``traditional'' codimension of \eqref{model-Hopf-for-def} is equal to the fractal codimension of \eqref{model-Hopf-for-def}, introduced in Definition \ref{def-intrin-codimen}. More precisely,
\begin{lemma}
\label{lemma-trad-frac}
Consider a smooth slow-fast Li\'{e}nard family \eqref{model-Hopf-for-def} with a finite codimension $j+1\ge 1$ (resp. infinite codimension), in the sense of \cite{DRbirth}. Then the fractal codimension of \eqref{model-Hopf-for-def} for $\lambda=\lambda_0$ is finite and equal to $j+1$ (resp. infinite).
\end{lemma}
Lemma \ref{lemma-trad-frac} will be proved in Section \ref{subsection-lemma-proof}.
\smallskip

From Theorem \ref{theorem-Lienard} it follows that the results from Remark \ref{remark-conjecture} can be stated not only for a classical Li\'{e}nard form \eqref{model-Hopf-for-def} but also for any slow-fast Hopf point using the more general notion of fractal codimension.

\begin{theorem}
\label{theorem-Lienard}
Consider a smooth slow-fast family $X_{\epsilon,\lambda}=X_{0,\lambda}+\epsilon Q_\lambda+O(\epsilon^2)$ that has a slow-fast Hopf point $p$ at $\lambda_0$ and satisfies Assumption C'. Then the following statements are true. 
\begin{enumerate}
    \item If the fractal codimension of $p$ is equal to $1$, then $\cycl(X_{\epsilon,\lambda},p)\le 1$. The limit cycle, if it exists, is hyperbolic and attracting (resp. repelling) when the slow divergence integral $I$ from Assumption C' is positive (resp. negative).
    \item If $p$ has finite fractal codimension $j+1\ge 1$ and if there exists a local chart on $M$ around $p$ in which, upon multiplication by a smooth positive function, $X_{\epsilon,\lambda}$ is given in a normal form \eqref{normal form 1} where $h=0$ and $f$ also may depend on $\epsilon$ (i.e., we have a Li\'{e}nard normal form), then $\cycl(X_{\epsilon,\lambda},p)$ is finite and bounded by $j+1$.
    \item If $X_{\epsilon,\lambda}$ is analytic on an analytic surface $M$, then $\cycl(X_{\epsilon,\lambda},p)$ is finite. Moreover, if $p$ has finite fractal codimension $j+1\ge 1$, then $\cycl(X_{\epsilon,\lambda},p)\le j+1$.
\end{enumerate}
\end{theorem}
Theorem \ref{theorem-Lienard} will be proved in Section \ref{section-proof-Lienard}. We use Lemma \ref{lemma-trad-frac} in the proof of Theorem \ref{theorem-Lienard}.2 and Theorem \ref{theorem-Lienard}.3.

\section{Proof of the main results}
\label{section-proofs}

\subsection{Proof of Theorem \ref{theorem-1}}
\label{proof-theorem1}

Suppose that $\tilde X_{\epsilon}$ defined in \eqref{normal form 2}
satisfies Assumptions A–C from Section \ref{subsectionNintrinsic}.  Let $y_0\in ]0,y^*[$ be arbitrary and fixed ($y^*$ comes from Assumption C). Suppose that $\tilde X_{\epsilon}$ has finite codimension $j+1\ge 1$ (Definition \ref{def-codim}), i.e. \[\tilde{g}(x)+\tilde{g}(-x)=\alpha x^{2j}+O(x^{2j+2}), \ \alpha\ne 0,\]
where $g(x,0)=g_m x^m+x^{m+1}\tilde{g}(x)$, with a smooth function $\tilde g$ and $g_m=\pm 1$. The integrand in \eqref{SDI-contact} can be written as
\begin{align}\label{proof-Th3.1.1}
    \frac{1}{g(x,0)}(nx^{n-1})^2&=\frac{n^2x^{2n-2-m}}{g_m+P_{e}(x)+\frac{\alpha}{2}x^{2j+1}(1+O(x))}\nonumber\\
    &=\frac{n^2x^{2n-2-m}}{g_m+P_{e}(x)}-\frac{n^2\alpha}{2}x^{2n-m+2j-1}(1+O(x))
\end{align}
where $P_e$ is an even polynomial of degree less than or equal to $2j$ with $P_e(0)=0$ ($P_e=0$ for $j=0$) and $O(x)$-functions are smooth. In the last step in \eqref{proof-Th3.1.1} we used $g_m^2=1$. Since $m$ is odd (Assumption A) and $P_e$ is even, the first term in \eqref{proof-Th3.1.1} is an odd function. Let's recall that $I(y,y)$ (denoted by $I(y)$) is the slow divergence integral defined in \eqref{SDI-contact}. We have 
\begin{align}\label{proof-Th3.1.2}
I(y)&=-\int_{-y^{1/n}}^{ y^{1/n}}\frac{1}{g(x,0)}(nx^{n-1})^2dx=\frac{n^2\alpha}{2}\int_{-y^{1/n}}^{ y^{1/n}}x^{2n-m+2j-1}(1+O(x))dx\nonumber\\
&=\frac{n^2\alpha}{2n-m+2j}y^{\frac{2n-m+2j}{n}}(1+o(1))
\end{align}
where $o(1)$ tends to zero as $y\to 0$. From \eqref{proof-Th3.1.2} it follows that $I$ is positive (resp. negative) if $\alpha>0$ (resp. $\alpha<0$).
Consider the sequence $(y_k)_{k\ge 0}$, starting with $y_0$, defined recursively by $I(y_{k+1},y_{k})=0$ or $I(y_{k},y_{k+1})=0$, $k\ge 0$, depending on the sign of $I$ and $g_m$ (see Section \ref{subsectionNintrinsic}). First we show that 
\begin{equation}\label{proof-Th3.1.3}
y_k-y_{k+1}\simeq y_k^\frac{n+2j+1}{n}, \ k\to \infty.
\end{equation}
We will prove \eqref{proof-Th3.1.3} when $I(y)>0$ on $]0,y^*[$ and $g_m=-1$, or $I(y)<0$ on $]0,y^*[$ and $g_m=+1$. The sequence $(y_k)_{k\ge 0}$ is then defined using $I(y_k,y_{k+1})=0$. The other case ($I>0$ and $g_m=+1$, or $I<0$ and $g_m=-1$) can be treated in a similar fashion. We have 
\begin{align}\label{proof-Th3.1.4}
I(y_k)&=-\int_{-y_k^{1/n}}^{ y_k^{1/n}}\frac{1}{g(x,0)}(nx^{n-1})^2dx      \nonumber\\
&=-\int_{-y_k^{1/n}}^{ y_k^{1/n}}\frac{1}{g(x,0)}(nx^{n-1})^2dx-I(y_k,y_{k+1})\nonumber\\
&=-\int_{y_{k+1}^{1/n}}^{ y_{k}^{1/n}}\frac{1}{g(x,0)}(nx^{n-1})^2dx\nonumber\\
&=-\frac{n^2}{g_m}\int_{y_{k+1}^{1/n}}^{ y_{k}^{1/n}}x^{2n-2-m}(1+O(x))dx\nonumber\\
&=-\frac{n}{g_m}\int_{y_{k+1}}^{ y_{k}}u^{\frac{n-m-1}{n}}(1+O(u^{\frac{1}{n}}))du,
\end{align}
for all $k$ sufficiently large. In the second step we used the fact that $I(y_k,y_{k+1})=0$, in the third step we used the definition of $I(\cdot,\cdot)$, given in \eqref{SDI-contact}, and in the last step we used the substitution  $u=x^n$. Using \eqref{proof-Th3.1.2} and \eqref{proof-Th3.1.4} we get
\begin{equation}\label{proof-Th3.1.5}
    \int_{y_{k+1}}^{ y_{k}}u^{\frac{n-m-1}{n}}(1+O(u^{\frac{1}{n}}))du\simeq y_k^{\frac{2n-m+2j}{n}}, \ k\to \infty.
\end{equation}
Assumption B implies that $\frac{n-m-1}{n}>-1$. If $\frac{n-m-1}{n}\ge 0$ (resp. $\frac{n-m-1}{n}<0$), then, on $[y_{k+1},y_k]$, the integrand in \eqref{proof-Th3.1.5} can be bounded from above by $y_k^{\frac{n-m-1}{n}}$ (resp. $y_{k+1}^{\frac{n-m-1}{n}}$) and from below by $y_{k+1}^{\frac{n-m-1}{n}}$ (resp. $y_k^{\frac{n-m-1}{n}}$), up to multiplication by a positive constant. Using this and \eqref{proof-Th3.1.5} we obtain
\begin{align}\label{proof-Th3.1.6}
    \frac{1}{\rho}\le \frac{y_k-y_{k+1}}{y_k^\frac{n+2j+1}{n}}\le \rho (\frac{y_k}{y_{k+1}})^{\frac{n-m-1}{n}}\nonumber\\
    \Big(\text{resp. } \frac{1}{\rho}(\frac{y_k}{y_{k+1}})^{\frac{n-m-1}{n}}\le \frac{y_k-y_{k+1}}{y_k^\frac{n+2j+1}{n}}\le \rho \Big)
\end{align}
 for some $\rho>0$ and all $k$ large enough. Now, it suffices to prove that
 \begin{align}\label{proof-Th3.1.7}
    \lim_{k\to\infty} \frac{y_k}{y_{k+1}}=1.
\end{align}
Then \eqref{proof-Th3.1.6} and \eqref{proof-Th3.1.7} imply \eqref{proof-Th3.1.3}. Let's prove \eqref{proof-Th3.1.7}. We make in the integral in \eqref{proof-Th3.1.5} the change of variable $u=y_k\tilde u$. We get 
\begin{equation}
    \int_{\frac{y_{k+1}}{y_k}}^{ 1}\tilde u^{\frac{n-m-1}{n}}(1+O((y_k\tilde u)^{\frac{1}{n}}))d\tilde u\simeq y_k^{\frac{2j+1}{n}}, \ k\to \infty,\nonumber
\end{equation}
and, consequently,
\begin{equation}
    \lim_{k\to\infty}\int_{\frac{y_{k+1}}{y_k}}^{ 1}\tilde u^{\frac{n-m-1}{n}}(1+O((y_k\tilde u)^{\frac{1}{n}}))d\tilde u=0.\nonumber
\end{equation}
Since the term $(1+O((y_k\tilde u)^{\frac{1}{n}}))$ is bounded from below by a positive constant, we have 
\begin{equation}
    \lim_{k\to\infty}\int_{\frac{y_{k+1}}{y_k}}^{ 1}\tilde u^{\frac{n-m-1}{n}}d\tilde u=0.\nonumber
\end{equation}
We can directly compute the above integral and see that \eqref{proof-Th3.1.7} is true.
\smallskip

Since $\eta:=\frac{n+2j+1}{n}>1$ in \eqref{proof-Th3.1.3}, \cite[Theorem 1]{EZZ} and \cite[Remark 1]{EZZ} imply that the set $(y_k)_{k\ge 0}$ is Minkowski nondegenerate, 
\[\dim_B(y_k)_{k\ge 0}=1-\frac{1}{\eta}=\frac{2j+1}{n+2j+1},\]
\[y_k\simeq k^{-\frac{1}{\eta-1}}, \ k\to \infty,\] and these results do not depend on the choice of $y_0$.
Since the map $y_k\mapsto (0,y_k)$ is bi-Lipschitz, we have that the fractal sequence $\mathcal S$ is
Minkowski nondegenerate and that $\dim_B\mathcal S=\dim_B(y_k)_{k\ge 0}$ (see Section \ref{section-Minkowski-def}). From formula (3.6.3) and Example 5.5.19 in \cite{Goran}, and from the above asymptotic for $y_k$, it follows that the chirp $\mathcal C_{\mathcal S}$ is Minkowski nondegenerate and 
\[\dim_B \mathcal C_{\mathcal S}=\frac{n+4j+1}{n+2j+1}.\]
\smallskip

When the codimension of the contact point is infinite, then instead of \eqref{proof-Th3.1.2} we have that for every $\tilde\eta>0$: $I(y)=O(y^{\tilde\eta})$, $y\to 0$. Then, using similar steps like in the finite codimension case, it can be seen that for any $\tilde\eta>1$: $y_k-y_{k+1}=O(y_k^{\tilde\eta})$, $k\to\infty$. \cite[Theorem 6]{EZZ} now implies that $\dim_B \mathcal{S}=1$ (the Minkowski dimension does not depend on the chosen initial condition $y_0$). This completes the proof of Theorem \ref{theorem-1}. 

\subsection{Proof of Theorem \ref{theorem-Main}}
\label{section-proof-Main}
Let a smooth slow-fast family $X_{\epsilon,\lambda}=X_{0,\lambda}+\epsilon Q_\lambda+O(\epsilon^2)$ have a contact point $p$ at $\lambda_0$ and satisfy Assumptions A, B and C'. Since the notions
of (upper, lower) Minkowski dimension and Minkowski nondegeneracy are defined in a coordinate-free way in Sections \ref{section-Minkowski-def} and \ref{subsection-riemann}, it suffices to prove Theorem  \ref{theorem-Main} in any local coordinates near $p$ for smooth equivalence, at level $\lambda=\lambda_0$. We may choose the normal form coordinates $(x,y)$ in a neighborhood of $p$ in which $X_{\epsilon,\lambda_0}$ is given by \eqref{normal form 2}, up to smooth equivalence (this means smooth coordinate change and a multiplication by a smooth positive function). We use the standard Euclidean metric. For the sake of readability we recall here the slow-fast family \eqref{normal form 2}:
\begin{equation}
    \widetilde{X}_\epsilon: \ \begin{vf}
        \dot{x} &=& y-x^{n} \nonumber \\
        \dot{y} &=&\epsilon\left( g(x,\epsilon)+\left(y-x^n\right)h(x,y,\epsilon)\right),\nonumber
    \end{vf}
\end{equation}
where $g$ and $h$ are smooth functions and $g(x,0)=g_mx^m+O(x^{m+1})$, with $g_m=\pm 1$. Following Theorem \ref{theorem-1}, if we choose $\sigma=\sigma_1\cup\{(0,0)\}\cup \sigma_2\subset\{x=0\}$ near the contact point $(x,y)=(0,0)$ (see Fig. \ref{fig-sequence}), then for any fractal sequence $\mathcal S$ (i.e. for any initial point $p_0=(0,y_0)\in\sigma_1$) we have 
\[\dim_B\mathcal S\in \left\{\frac{2j+1}{n+2j+1} \ | \ j\in\mathbb N_0 \right\} \text{ or } \dim_B\mathcal S=1. \]
Moreover, if $\dim_B\mathcal S\ne 1$, then $\mathcal S$ is Minkowski nondegenerate.

Now, let us show that the Minkowski dimension and Minkowski nondegeneracy for fractal sequences remain invariant if we replace $\sigma$ with a new smooth regular section $\sigma'=\sigma_1'\cup\{(0,0)\}\cup \sigma_2'$, transverse to horizontal lines. Since $\sigma'$ can be parameterized in a smooth way by the variable $y$, $\sigma'$ is the graph of a smooth function $x=\tilde{h}(y)$, and the map $p_k=(0,y_k)\in\mathcal S\mapsto p_k'=(\tilde{h}(y_k),y_k)\in \mathcal S'$ is therefore bi-Lipschitz. Now, using Section \ref{section-Minkowski-def}, it follows that the Minkowski dimension and Minkowski nondegeneracy for fractal sequences generated by $\tilde X_{\epsilon}$  do not depend on the choice of the section $\sigma'$ and the initial point on it. 

Suppose now that $\dim_B\mathcal S\ne 1$ for some fractal sequence $\mathcal S$ generated by $\tilde X_{\epsilon}$. Then the same is true for any fractal sequence generated by $\tilde X_{\epsilon}$ using the above result. Theorem \ref{theorem-1} and Remark \ref{remarkSS'} imply that chirp produced by any fractal sequence is Minkowski nondegenerate and has the Minkowski dimension which does not change as we vary fractal sequences. This completes the proof of Theorem \ref{theorem-Main}.

\subsection{Proof of Lemma \ref{lemma-trad-frac}}
\label{subsection-lemma-proof}
We consider a smooth Li\'{e}nard system \eqref{model-Hopf-for-def} and write \[\Psi_{\epsilon,\lambda}(x)=x\sqrt{1-xh_1(x,\epsilon,\lambda)}\] where $h_1$ is defined in \eqref{model-Hopf-for-def}. After the change of coordinates $\bar x=\Psi_{\epsilon,\lambda}(x)$ and division by $\Psi_{\epsilon,\lambda}'(x)>0$, the system \eqref{model-Hopf-for-def} changes near $(x,y)=(0,0)$ into
\begin{equation}
\label{model-Hopf-for-def-2}
    \begin{vf}
        \dot{\bar x} &=& y-\bar x^2 \\
        \dot{y} &=&\epsilon g(\bar x,\epsilon,\lambda),
    \end{vf}
\end{equation}
where $g$ is a smooth function and 
\begin{equation}
    \label{Liena-formula-1}
    g(\bar x,0,\lambda_0)=-\Psi^{-1}(\bar x)(\Psi^{-1})'(\bar x)
\end{equation}
with $\Psi( x):=\Psi_{0,\lambda_0}( x)$. We can write $g(\bar x,0,\lambda_0)=-\bar x+\bar x^2\tilde g(\bar x)$ where $\tilde g$ is smooth. Note that for $\lambda=\lambda_0$ system \eqref{model-Hopf-for-def-2} is of type \eqref{normal form 2} with $(n,m)=(2,1)$. If \eqref{model-Hopf-for-def} has a finite codimension $j+1$, in the sense of \cite{DRbirth}, then we will prove that 
\begin{equation}\label{tilda-g-lemma}\tilde{g}(\bar x)+\tilde{g}(-\bar x)=\alpha \bar x^{2j}+O(\bar x^{2j+2}), 
\end{equation}
for some $\alpha\ne 0$. This implies (see Definition \ref{def-codim}) that the fractal codimension of the normal form \eqref{model-Hopf-for-def-2} is equal to $j+1$. As the notion of fractal codimension is defined in a coordinate-free way (Section \ref{subsection-int-def-codim}), we have that the fractal codimension of \eqref{model-Hopf-for-def} is equal to $j+1$. If \eqref{model-Hopf-for-def} has infinite codimension, in the sense of \cite{DRbirth}, then we will show that $j\in\mathbb N_0$ with the property \eqref{tilda-g-lemma} does not exist. Then the fractal codimension of \eqref{model-Hopf-for-def-2} is infinite (Definition \ref{def-codim}). Again, using the intrinsic nature of the notion of fractal codimension, we conclude that the fractal codimension of \eqref{model-Hopf-for-def} is infinite.
\smallskip

If \eqref{model-Hopf-for-def} has a finite codimension $j+1$, then \eqref{h1even} and the definition of $\Psi$ imply that 
\begin{equation}
    \label{vazno-lemma-1}
    \Psi(x)= xP_e( x)+\alpha  x^{2j+2}(1+O( x))
\end{equation}
where $P_e$ is an even polynomial of degree less than or equal to $2j$, $P_e(0)=1$, and where $\alpha$ is a new nonzero constant. On the other hand, if \eqref{model-Hopf-for-def} has infinite codimension, then for any fixed $\tilde j\in\mathbb{N}_0$ we get 
\begin{equation}
    \label{vazno-lemma-2}
    \Psi(x)= x\tilde P_e( x)+O(x^{2\tilde j+3})
\end{equation}
where $\tilde P_e$ is an even polynomial of degree less than or equal to $2\tilde j$ and $\tilde P_e(0)=1$. We show that $\Psi^{-1}$ can be written as \eqref{vazno-lemma-1} in finite codimension case, for some new $P_e$ and $\alpha$ with the same property as in \eqref{vazno-lemma-1}, or as \eqref{vazno-lemma-2} in infinite codimension case, for some new $\tilde P_e$ with the same property as in \eqref{vazno-lemma-2}. To prove this, we apply The Lagrange Inversion Formula for $\Psi(x)$: for each positive integer $i$, the coefficient of $\bar x^i$ in the Taylor
expansion at $\bar x = 0$ of the inverse $\Psi^{-1}(\bar x)$ is equal to the coefficient of $ x^{i-1}$ in the Taylor
expansion at $x = 0$ of 
\[\left(\frac{x}{\Psi(x)}\right)^i,\]
divided by $i$.

First, we focus on the finite codimension case. For each fixed $i$, it can be easily seen, using \eqref{vazno-lemma-1}, that 
\begin{equation}
    \label{vazno-lemma-3}
    \left(\frac{x}{\Psi(x)}\right)^i=P_e(x)+\alpha x^{2j+1}(1+O(x)),
\end{equation}
for some new $P_e$ and $\alpha$ with the same property as in \eqref{vazno-lemma-1}.  It follows from \eqref{vazno-lemma-3} and The Lagrange Inversion Formula that the coefficient of $\bar x^i$ in $\Psi^{-1}(\bar x)$ is zero, for each $i=2,\dots,2j$, and nonzero for $i=2j+2$. Indeed, when $i$ is even, then the $(i-1)$st derivative $P_e^{(i-1)}(0)=0$ because $P_e$ in \eqref{vazno-lemma-3} is even. On the other hand, the $(i-1)$st derivative of the second term in \eqref{vazno-lemma-3} calculated at $x=0$ is zero for $i=2,\dots,2j$ and nonzero for $i=2j+2$ ($\alpha\ne 0$). This implies that the coefficient of $x^{i-1}$ in \eqref{vazno-lemma-3} is $0$ for each $i=2,\dots,2j$ and nonzero for $i=2j+2$. Thus, using The Lagrange Inversion Formula, we see that $\Psi^{-1}$ has an expansion of type \eqref{vazno-lemma-1} with new $P_e$ ($P_e(0)=1$) and nonzero $\alpha$. Then we get ($g(\bar x,0,\lambda_0)$ is given in \eqref{Liena-formula-1}): 
\begin{align*}
    g(\bar x,0,\lambda_0)&=-\Psi^{-1}(\bar x)(\Psi^{-1})'(\bar x)\\
    &=-\bar x\left(P_e( \bar x)+\alpha  \bar x^{2j+1}(1+O( \bar x))\right)\\
    & \qquad\qquad \left(P_e( \bar x)+\bar xP_e'( \bar x)+\alpha  (2j+2)\bar x^{2j+1}(1+O( \bar x))\right)\\
    &=-\bar x\left(P_E( \bar x)+\alpha(2j+3)  \bar x^{2j+1}(1+O( \bar x))\right)\\
    &=-\bar x+\bar x^2\left(P_o( \bar x)-\alpha(2j+3)  \bar x^{2j}(1+O( \bar x)) \right),
\end{align*}
where $P_E$ (resp. $P_o$) is an even (resp. odd) polynomial of degree less than or equal to $2 j$ (resp. $2j-1$) and $P_E(0)=1$. Thus, we have proved \eqref{tilda-g-lemma}.

In the infinite codimension case, we use \eqref{vazno-lemma-2} and obtain, for each fixed $i\in\mathbb N$ and $\tilde j\in\mathbb N_0$, 
\begin{equation}
    \label{vazno-lemma-4}
    \left(\frac{x}{\Psi(x)}\right)^i=\tilde P_e(x)+ O(x^{2\tilde j+2}),
\end{equation}
for some new $\tilde P_e$ with the same property as in \eqref{vazno-lemma-2}. The Lagrange Inversion Formula and \eqref{vazno-lemma-4} imply now that the coefficient of $\bar x^i$ in $\Psi^{-1}(\bar x)$ is zero, for each $i=2,\dots,2\tilde j+2$ (see the finite codimension case). Thus, for each fixed $\tilde j\in\mathbb N_0$, $\Psi^{-1}$ has an expansion of type \eqref{vazno-lemma-2}, for some new $\tilde P_e$ ($\tilde P_e(0)=1$). We have for any fixed $\tilde j\in\mathbb N_0$:
\begin{align*}
    g(\bar x,0,\lambda_0)&=-\Psi^{-1}(\bar x)(\Psi^{-1})'(\bar x)\\
    &=-\bar x\left(\tilde P_e( \bar x)+O(  \bar x^{2\tilde j+2})\right)
    \left(\tilde P_e( \bar x)+\bar x\tilde P_e'( \bar x)+O(  \bar x^{2\tilde j+2})\right)\\
    &=-\bar x\left(\tilde P_E( \bar x)+O(  \bar x^{2\tilde j+2})\right)\\
    &=-\bar x+\bar x^2\left(\tilde P_o( \bar x)+O(  \bar x^{2\tilde j+1}) \right),
\end{align*}
where $\tilde P_E$ (resp. $\tilde P_o$) is an even (resp. odd) polynomial of degree less than or equal to $2 \tilde j$ (resp. $2\tilde j-1$) and $\tilde P_E(0)=1$. This implies that \eqref{tilda-g-lemma} can never hold. We have proved Lemma \ref{lemma-trad-frac}.

\subsection{Proof of Theorem \ref{theorem-Lienard}}
\label{section-proof-Lienard}
Let $X_{\epsilon,\lambda}=X_{0,\lambda}+\epsilon Q_\lambda+O(\epsilon^2)$ be a smooth slow-fast family with a slow-fast Hopf point $p\in M$ at $\lambda_0$. Assume that $p$ satisfies Assumption C' (we know that $p$ satisfies Assumptions A and B). Since $p$ is a slow-fast Hopf point, we can use the normal form \eqref{normal form 1}, for $(\epsilon,\lambda)\sim (0,\lambda_0)$, where $f(0,\lambda_0)=\frac{\partial f}{\partial x}(0,\lambda_0)=0$, $\frac{\partial^2 f}{\partial x^2}(0,\lambda_0)\ne 0$, $g(0,0,\lambda_0)=0$ and $\text{sign} \frac{\partial g}{\partial x}(0,0,\lambda_0)=- 1$. 

\textit{Proof of Theorem \ref{theorem-Lienard}.1.} The contact point $p$ is of Morse type and, after translations in $x$ and $y$ and a rescaling in $(x,y)$, \eqref{normal form 1} changes into
\begin{equation}
\label{normal form 1'}
    \begin{vf}
        \dot{x} &=& y-f(x,\lambda)  \\
        \dot{y} &=&\epsilon\left( g(x,\epsilon,\lambda)+\left(y-f(x,\lambda)\right)h(x,y,\epsilon,\lambda)\right),
    \end{vf}
\end{equation}
for some new smooth functions $f$, $g$ and $h$ where $f(x,\lambda)=x^2+x^3f_1(x,\lambda)$ and $g$ has the property mentioned above (see \cite[Section 6.1]{DDR-book-SF}). Using the $\lambda$-dependent coordinate change $\bar x=x\sqrt{1+xf_1(x,\lambda)}$, the normal form \eqref{normal form 1'}, near $(x,y)=(0,0)$, can be written as 
\begin{equation}
\label{normal form 2'}
    \begin{vf}
        \dot{x} &=& y-x^2  \\
        \dot{y} &=&\epsilon\left( g(x,\epsilon,\lambda)+\left(y-x^2\right)h(x,y,\epsilon,\lambda)\right),
    \end{vf}
\end{equation}
up to multiplication by a smooth positive function, where $g$ and $h$ are new smooth functions, $g$ has the same properties and we write $x$ instead of $\bar x$. Using an $(\epsilon,\lambda)$-dependent rescaling in $(x,y,t)$ we can assume that $g(0,0,\lambda_0)=0$ and $ \frac{\partial g}{\partial x}(0,\epsilon,\lambda)=- 1$, for all $(\epsilon,\lambda)\sim (0,\lambda_0)$, in \eqref{normal form 2'}. System \eqref{normal form 2'} is a normal form for $C^\infty$-equivalence and it is valid in an $(\epsilon,\lambda)$-uniform neighborhood of $(x,y)=(0,0)$.

When $\lambda=\lambda_0$, the normal form \eqref{normal form 2'} is of type \eqref{normal form 2} with $(n,m)=(2,1)$. Since we suppose that the codimension of $p$ in $X_{\epsilon,\lambda}$ is equal to $1$, then Definition \ref{def-intrin-codimen} implies that $\dim_B\mathcal S=\frac{1}{3}$ for any fractal sequence $\mathcal S$ generated by the normal form \eqref{normal form 2'} (see also Theorem \ref{theorem-Main}). Now Theorem \ref{theorem-1} and Definition \ref{def-codim}  imply that $ \frac{\partial^2 g}{\partial x^2}(0,0,\lambda_0)\ne 0$ in \eqref{normal form 2'}. Moreover, from Remark \ref{remark-alpha} it follows that $ \frac{\partial^2 g}{\partial x^2}(0,0,\lambda_0)> 0$ (resp. $ \frac{\partial^2 g}{\partial x^2}(0,0,\lambda_0)<0$) if the slow divergence integral $I$ from Assumption C' is positive (resp. negative). Note that we use the fact that the notion of slow divergence integral is invariant under $C^\infty$-equivalences (see Section \ref{subsub-general}).
Following \cite[Section 4.4]{CriticalityHopf}, we have that the cyclicity of the origin in \eqref{normal form 2'} is bounded by one and the limit cycle, if it exists, is hyperbolic and stable when $ \frac{\partial^2 g}{\partial x^2}(0,0,\lambda_0)> 0$ and hyperbolic and unstable when $ \frac{\partial^2 g}{\partial x^2}(0,0,\lambda_0)<0$. This completes the proof of Theorem \ref{theorem-Lienard}.1.
\smallskip

\textit{Proof of Theorem \ref{theorem-Lienard}.2.} We suppose that $p$ has finite fractal codimension $j+1\ge 1$ with 
\[j=\frac{3\dim_B\mathcal S-1}{2(1-\dim_B\mathcal S)}\]
where $\mathcal S$ is any fractal sequence near $p$ (see Definition \ref{def-intrin-codimen}). We also suppose that there exists a local chart on $M$ around $p$ in which, upon multiplication by a smooth positive function, $X_{\epsilon,\lambda}$, with $(\epsilon,\lambda)\sim (0,\lambda_0)$, can be written as a Li\'{e}nard system
\begin{equation}
\label{normal form Lienar}
    \begin{vf}
        \dot{x} &=& y-f(x,\epsilon,\lambda)  \\
        \dot{y} &=&\epsilon  g(x,\epsilon,\lambda)
    \end{vf}
\end{equation}
where $p$ is given by $(x,y)=(0,0)$ and $f$ and $g$ are smooth functions. Since $p$ is a slow-fast Hopf point, we know that $f(0,0,\lambda_0)=\frac{\partial f}{\partial x}(0,0,\lambda_0)=0$, $\frac{\partial^2 f}{\partial x^2}(0,0,\lambda_0)\ne 0$, $g(0,0,\lambda_0)=0$ and $ \frac{\partial g}{\partial x}(0,0,\lambda_0)<0$. Following \cite{DRbirth}, the system \eqref{normal form Lienar} is smoothly equivalent to a ``classical'' Li\'{e}nard equation of type \eqref{model-Hopf-for-def}. For the sake of completeness, we include here a slightly different proof of this fact. Then, in the classical Li\'{e}nard setting, 
 we can use cyclicity results of \cite{DRbirth}. After an $(\epsilon,\lambda)$-dependent translation in $x$ and a rescaling in $(x,t)$, the system \eqref{normal form Lienar}, near $(x,y)=(0,0)$, becomes 
\begin{equation}
\label{normal form Lienar3}
    \begin{vf}
        \dot{x} &=& y-f(x,\epsilon,\lambda)  \\
        \dot{y} &=&\epsilon (-x+x^2\tilde{g}(x,\epsilon,\lambda)), 
    \end{vf}
\end{equation}
where $ f$ and $\tilde g$ are (new) smooth functions and $f$ has the same property as in \eqref{normal form Lienar}. The function $1-x\tilde{g}(x,\epsilon,\lambda)$ is strictly positive near $x=0$, and there is a strictly positive smooth function $\tilde G(x,\epsilon,\lambda)$, $\tilde G(0,\epsilon,\lambda)=1$, such that \[\frac{x^2}{2}\tilde G(x,\epsilon,\lambda)=\int_0^xs(1-s\tilde{g}(s,\epsilon,\lambda))ds.\]
If we differentiate this w.r.t. $x$, then we get 
\begin{equation}
    \label{diff-G}
    \tilde G(x)+\frac{x}{2}\tilde G'(x)=1-x\tilde{g}(x),
\end{equation}
after division by $x$. If we use the coordinate change $\bar x=x\sqrt{\tilde G(x,\epsilon,\lambda)}$, \eqref{diff-G} and multiplication by 
\[\frac{\sqrt{\tilde G(x,\epsilon,\lambda)}}{1-x\tilde{g}(x,\epsilon,\lambda)}>0,\]
then the system \eqref{normal form Lienar3}, near $(x,y)=(0,0)$, changes into
\begin{equation}
\label{normal form Lienar4}
    \begin{vf}
        \dot{ x} &=& y-f(x,\epsilon,\lambda)  \\
        \dot{y} &=&-\epsilon x, 
    \end{vf}
\end{equation}
where we denote $\bar x$ again by $x$ and $f$ is a new smooth function having the same property as in \eqref{normal form Lienar}.  Finally, after an $(\epsilon,\lambda)$-dependent translation in $(x,y)$ and a linear rescaling, we can bring \eqref{normal form Lienar4} into a smooth system
\begin{equation}
\label{normal form Lienar5}
    \begin{vf}
        \dot{ x} &=& y-(x^2+O(x^3))  \\
        \dot{y} &=&\epsilon (b(\epsilon,\lambda)-x), 
    \end{vf}
\end{equation}
with $b(0,\lambda_0)=0$.
Clearly, \eqref{normal form Lienar5} is in a classical Li\'{e}nard setting \eqref{model-Hopf-for-def}.

Note that the fractal codimension of \eqref{normal form Lienar5} is $j+1$ because of its intrinsic nature. Now, using Lemma \ref{lemma-trad-frac}, we know that the codimension of \eqref{normal form Lienar5} (in the sense of \cite{DRbirth}) is also equal to $j+1$.
Following \cite[Theorem 1.5]{DRbirth}, \eqref{normal form Lienar5} has a finite cyclicity at $(x, y)=(0,0)$. Moreover, the cyclicity is bounded by $j+1$ up to the conjecture from Remark \ref{remark-conjecture}. The conjecture has been proved in \cite{Li-conjecture-solved}. Since we have proved that there are local coordinates $( x,y)$ in which $p=(0,0)$ and $X_{\epsilon,\lambda}$, with $(\epsilon,\lambda)\sim (0,\lambda_0)$, is equal to \eqref{normal form Lienar5}, up to smooth equivalence, we have that $\cycl(X_{\epsilon,\lambda},p)$ is finite and bounded by $\cdim(p)=j+1$. This finishes the proof of Theorem \ref{theorem-Lienard}.2.
\smallskip

\textit{Proof of Theorem \ref{theorem-Lienard}.3.} If $X_{\epsilon,\lambda}$ is an analytic slow-fast family with a slow-fast Hopf point $p$, then \cite[Theorem 4]{RHAnalytic} implies that the family $X_{\epsilon,\lambda}$, near $p$ and for $(\epsilon,\lambda)\sim (0,\lambda_0)$, is given by an analytic Li\'{e}nard family of type \eqref{normal form Lienar}, up to analytic equivalence. Following \cite[Theorem 1.2]{DRbirth}, the analytic Li\'{e}nard family with a slow-fast Hopf point has a finite cyclicity, i.e. $\cycl(X_{\epsilon,\lambda},p)$ is finite. If $p$ has finite fractal codimension $\cdim(p)=j+1\ge 1$, then $\cycl(X_{\epsilon,\lambda},p)\le j+1$ (see Theorem \ref{theorem-Lienard}.2).

\section{Minkowski dimension of subsets of Riemannian manifolds}
\label{subsection-riemann}

Since we are working with fractal objects defined on Riemannian manifolds, in this subsection we show that the notion of Minkowski dimension can be intrinsically defined in this context.
More specifically, the Minkowski dimension and Minkowski nondegeneracy will be independent of the choice of the Riemannian metric.

Assume from now on that $(M,g)$ is a compact Riemann manifold, i.e., a compact real smooth manifold $M$ where $g_p$ is a positive-definite inner product on the tangent space $T_pM$ at each point $p\in M$.

The Riemann metric tensor $g_p$, i.e., the family $\{g_p:p\in M\}$ then induces a structure of a metric space on the smooth manifold $M$ in the standard way via the induced {\em Riemann distance function} $d_g$ where $d_g(x,y)$ is the infimum of lengths of all admissible\footnote{A piecewise continuously differentiable curve from $x$ to $y$ with nonzero velocity, in fact somewhat more general but not important for our discussion; see, e.g., \cite{lee-book} for more details.
The length of the curve is defined in a standard way by integrating the norm (induced by $g_p$) of the velocity along the curve.} curves from $x$ to $y$ 

It is well known (see, e.g., \cite[Chapter 2]{lee-book}) that the topology induced by $d_g$ on $M$ coincides with the manifold topology of $M$.
Moreover, since $M$ is compact, then any two Riemann distance functions are uniformly equivalent, i.e., if $d_g$ and $d_h$ are distance functions induced by Riemann metric tensors $g_p$ and $h_p$, respectively, then there exist positive constants $A$ and $B$ such that
\begin{equation}
	\label{rim-eq}
	A\,d_h(x,y)\leq d_g(x,y)\leq B\,d_h(x,y)\quad\textrm{for all}\quad x,y\in M.
\end{equation}

Consider now the compact Riemann manifold $(M,g)$ as a metric measure space $(M,d_g,\mathcal{H}_g^s)$ where $\mathcal{H}_g^s$ is the $s$-dimensional Hausdorff measure on $M$ induced by the metric $d_g$.
More precisely, let $F\subseteq M$ and
\begin{equation}
	\label{hausdorff-dm}
	\mathcal{H}_{\delta,g}^s(F):=\inf\left\{\sum_{i=1}^{\infty}(\mathrm{diam}_g(U_i))^s: \{U_i\}_{i=1}^{\infty}\ \textrm{is a $\delta$-cover of}\ F\right\},
\end{equation}
where $s\geq 0$, $\delta >0$, $\mathrm{diam}_g$ denotes the diameter of the set in the metric $d_g$ and all of $U_i$ have diameter at most $\delta$.
Then one defines the \emph{$s$-dimensional Hausdorff measure of $F$} as
\begin{equation}
	\label{hausdorff-m}
	\mathcal{H}_{g}^s(F):=\lim_{\delta\to 0^+}\mathcal{H}_{\delta,g}^s(F).
\end{equation}
The above limit exists for any subset of $M$ and defines a measure on the $\sigma$-algebra of Borel subsets of $M$.\footnote{The construction is a special case of Carath\'eodory's method, see, e.g., \cite{Falconer,Rog-H-measures,Mattila}.}

Note that the Hausdorff measure behaves nicely under Lipschitz mappings in the sense that if $F\subseteq M$ and $f\colon F\to M$ is Lipschitz then $\mathcal{H}^s_g(f(F))\leq {C^s}\mathcal{H}^s_g(F)$ where $C$ is the Lipschitz constant of $f$.\footnote{
Strictly speaking, the left-hand side of the inequality needs to be understood in the sense of Hausdorff outer measure since $f(F)$ need not be Borel in general if $f$ is just Lipschitz. However, if $f$ is bijective, i.e., bi-Lipschitz, there is no problem.}
This, in turn, implies that changing the Riemann metric tensor $g_p$ will result in equivalent Hausdorff measures.
Indeed, as already discussed, on compact Riemann manifolds any two Riemann distance functions, say, $d_g$ and $d_h$ are uniformly equivalent, i.e., \eqref{rim-eq} holds, then the identity map from $(M,d_g)$ to $(M,d_h)$ is bi-Lipschitz and hence,
for any Borel subset $F$ of $M$ we have
\begin{equation}
	\label{haus-m-eq}
	A^s\,\mathcal{H}_h^s(F)\leq\mathcal{H}_g^s(F)\leq B^s\, \mathcal{H}_h^s(F).
\end{equation}

Since $M$ is a smooth manifold the only Hausdorff measure on $M$ of interest, that is, the one which is nontrivial, is in the case when $s$ equals to $N$, the topological dimension of $M$.
Namely, $M$ is locally bi-Lipshitz equivalent to $\mathbb{R}^N$ and hence, we also have local equivalence of the corresponding $s$-dimensional Hausdorff measures.

Finally, we now define the intrinsic Minkowski dimension for subsets of compact Riemannian manifolds.
\begin{definition}
	\label{riemann-mink}
	Let $(M,g)$ be a compact Riemann manifold of topological dimesion $N$, $F\subseteq M$ and $r\geq0$. We define the upper $r$-dimensional $g$-Minkowski content of $f$ by
	\begin{equation}
		\label{riemann-m-eq}
		\overline{\mathcal{M}}^r_g(F):=\limsup_{\varepsilon\to0^+}\frac{\mathcal{H}^N_g(F_{g,\varepsilon})}{\varepsilon^{N-r}}
	\end{equation}
	where $F_{g,\varepsilon}$ denotes the $\varepsilon$ neighborhood of $F$ in the Riemann distance $d_g$.
	
	Furthermore, the upper intrinsic Minkowski dimension of $F$ is then defined as
	\begin{equation}
		\label{upper-intr}
		\overline{\dim}_{B,M}:=\inf\{r:\overline{\mathcal{M}}^r_g(F)=0\}=\sup\{r:\overline{\mathcal{M}}^r_g(F)=+\infty\}.
	\end{equation}
		
\end{definition}
	
One now defines analogously as in the classical setting the notions of lower intrinsic Minkowski content and intrinsic Minkowski dimension.
In light of the just discussed equivalence of Riemann distance functions and the associated Hausdorff measures it is rather straightforward to show that the notions of (upper, lower) intrinsic Minkowski dimension do not depend on the choice of the Riemann metric tensor and this is also true for the notion of intrinsic Minkowski nondegeneracy (defined analogously as in the classical setting).

Moreover, if one wants to determine the intrinsic Minkowski dimension (and intrinsic Minkowski nondegeneracy) one can work in charts since the chart maps are locally bi-Lipschitz.
More precisely, if the set $A\subseteq M$ is ``sufficiently local'', i.e., such that there exists a chart $(U,\varphi)$ such that $A\subseteq U$ and $\varphi\colon A\to\varphi(A)$ is bi-Lipschitz then the intrinsic (upper, lower) Minkowski dimension is equal to the classical (upper, lower) Minkowski dimension of $\varphi(A)$ and also $A$ is intrinsically Minkowski nondegenerate if and only if $\varphi(A)$ is clasically Minkowski nondegenerate.\footnote{
If this is not the case one can always partition $A$ into smaller parts such that each part is ``sufficiently local'' and then can use the \emph{finite stability property} of the upper Minkowski dimension.
However, note the finite stability property does not hold for the lower Minkowski dimension which makes this situation more complicated and possibly interesting for further studies.}
Conveniently, in applications in this paper one is always able to choose a set $A$ which is ``sufficiently local''.
Namely, a sufficiently small neighborhood of the contact point $p$ will contain both sets on which we apply the theory, i.e., the fractal sequence $\mathcal{S}$ and the associated fractal chirp $\mathcal{C}_{\mathcal{S}}$.

\section{Applications}
\label{section-applications}
In this section we numerically compute the Minkowski dimension of fractal sequences for some planar contact points. Then we can read their fractal codimension from the Minkowski dimension (Definition \ref{def-intrin-codimen}). We will show that what was proven in the previous sections of the paper corresponds to the numerical data.
We give three examples. In Sections \ref{section-examp1}-\ref{section-examp2}, we deal with slow-fast Hopf points, inside classical Li\'{e}nard equations (Section \ref{section-examp1}) and a two-stroke oscillator discussed in \cite{Martin} (Section \ref{section-examp2}). In Section \ref{section-examp4} we have a contact point of type \eqref{normal form 2}, with a general contact order and singularity order.

\smallskip

For all of the examples shown here Assumptions A, B and C (or C') are fulfilled. To get these numerical results, the programming tool of Matlab was used. The code is open source and; see \url{https://github.com/AnsfriedJanssens/Code-fractal-dimensions-/tree/main}. To numerically calculate the Minko\-wski dimension of fractal sequences in two dimensional slow-fast systems, we use the following result (see \cite{BoxVlatko}):
\begin{proposition}\label{prop-3-meth} Let $\mathcal S = ( y_k )_{k \ge 0}$ be a decreasing sequence which tends to $0$ with an initial point
$y_0$ such that $y_k-y_{k+1} \simeq y_k^{\eta}$ as $k \rightarrow \infty$, for $\eta > 1$. Then $\dim_B\mathcal S=\frac{\eta-1}{\eta}$ and we have that
\begin{equation}\label{meth-1}
    \dim_B\mathcal S=\lim_{k \to \infty}  \frac{\ln k}{-\ln(y_k-y_{k+1})},
\end{equation}
\begin{equation}\label{meth-2}
    \dim_B\mathcal S=\lim_{k \to \infty} \frac{1}{1-\frac{\ln y_k}{\ln k}}
\end{equation}
and
\begin{equation}\label{meth-3}
    \dim_B\mathcal S=\lim_{k \to \infty}\left( 1- \frac{\ln(k(y_k-y_{k+1})+y_k)}{\ln(\frac{y_k-y_{k+1}}{2})}\right).
\end{equation}
\end{proposition}
\begin{remark}
The formulas given in \eqref{meth-1}, \eqref{meth-2} and \eqref{meth-3} are motivated by \cite[Section 3.4]{Tricot}. In \eqref{meth-1} we have a formula of Cahen-type. The formula given in \eqref{meth-2} is related to the so-called Borel
rarefaction index of $\mathcal S$ while \eqref{meth-3} is based on the decomposition of $\mathcal S$ into tail and nucleus (see also \cite{BoxVlatko}). 
\end{remark}
Roughly speaking, when we expect higher density of a fractal sequence $\mathcal S$ from Proposition \ref{prop-3-meth} (i.e., larger Minkowski dimension $\dim_B\mathcal S$), then the sequence given in \eqref{meth-3} produces the smallest error for a large fixed $k$. If $\dim_B\mathcal S$ is ``closer'' to zero, then one should use the sequence in \eqref{meth-1} or \eqref{meth-2} to estimate $\dim_B\mathcal S$. For more details see \cite[Figure 4]{BoxVlatko}. 
When showing the results only the method that works the best/fastest is shown.\\

In the examples we parameterize vertical sections with the variable $y$ in the $(x,y)$-phase space. We will deal with a decreasing or increasing sequence $(y_k)_{k \ge 0}$ that converges to $a\in\mathbb R$. Note that we need to have $\lim_{k\to \infty}y_k=0$ with $y_k>0$ $\forall k=0,1,2,\dots $ (Proposition \ref{prop-3-meth}). This gets solved by replacing $(y_k)_{k \ge 0}$ with $( |y_k-a| )_{k \ge 0}$ which is now decreasing and tends to zero. This bi-Lipschitz transformation doesn't change the Minkowski dimension (see Section \ref{section-Minkowski-def}). The starting point $y_0$ didn't influence the result of the methods. 

\subsection{Classical Li\'{e}nard equations}\label{section-examp1}
Let's first discuss the Minkowski dimension for a classical Li\'{e}nard equation
 \begin{equation}
\label{Application 1 system}
    \begin{vf}
        \dot{x} &=& y-\left(x^2+\alpha_{2j+3}x^{2j+3}\right) \\
        \dot{y} &=&-\epsilon x,
    \end{vf}
\end{equation}
where $j\in\mathbb N_0$, $\alpha_{2j+3} \neq 0$ and $\epsilon\ge 0$ is the singular perturbation parameter. The critical curve of \eqref{Application 1 system}, with $\epsilon=0$, is given by $\{y=F(x)\}$ where $F(x)=x^2+\alpha_{2j+3}x^{2j+3}$. System \eqref{Application 1 system} has a slow-fast Hopf point at $(x,y)=(0,0)$ (see \eqref{model-Hopf-for-def}). We define a vertical section $\sigma\subset \{x=0\}$, near $(x,y)=(0,0)$, parameterized by the variable $y$. If $(0,y), (0,\tilde y)\in\sigma$ and $y,\tilde y>0$, then the slow divergence integral along $[\alpha(y),\omega(\tilde y)]$ is given by 
\begin{equation}\label{Application 1 integral}
I(y,\tilde y)=\int^{\omega(\tilde y)}_{\alpha(y)}\frac{(F'(x))^2}{x}dx\nonumber
\end{equation}
where $x=\alpha(y)<0$ (resp. $x=\omega(\tilde y)>0$) is the $x$-component of the $\alpha$-limit point (resp. $\omega$-limit point) of the orbit of \eqref{Application 1 system}, with $\epsilon=0$, passing through $(0,y)$ (resp. $(0,\tilde y)$). Then, for a small $y_0\in \sigma$ ($y_0>0$), we generate a sequence $(y_k)_{k \ge 0}$ that converges to $0$ using either $I(y_k,y_{k+1})=0$ or $I(y_{k+1},y_k)=0$, depending on the sign of $\alpha_{2j+3}$ (see Section \ref{section-statement}).
Table \ref{table classical} gives for certain parameters the expected Minkowski dimensions and the results of the Matlab code for those parameters. Notice that only the value of the parameter $j$ decides the Minkowski dimensions, as our theory predicted. The parameter $\alpha_{2j+3}$ doesn't influence the Minkowski  dimension, it only takes more iterations for the sequences of the used methods to reach the Minkowski dimension for big $\alpha_{2j+3}$. Indeed, the system \eqref{Application 1 system} has finite codimension $j+1$ in the sense of \cite{DRbirth}. Lemma \ref{lemma-trad-frac} implies that the fractal codimension of \eqref{Application 1 system} is finite and equal to $j+1$. We obtain the same fractal codimension from the estimated Minkowski dimension (the last column of Table \ref{table classical}) if we use Definition \ref{def-intrin-codimen}.  The starting point for $y_0$ needs to be close to the contact point $y=0$. For higher values of $j$, we took the starting point $y_0=0.3$ further away from the contact point because it is more time efficient for the numerical calculations in the algorithm. 
\begin{table}
 \centering
\begin{tabular}{|c|c|c|c|c|c|c|}
\hline
\textbf{\# Iterations} &  \textbf{$y_0$} & $j$ &  \textbf{$\alpha_{2j+3}$} & \textbf{Theoretical Value} & \textbf{Results}\\
\hline
100 & 0.001 & 0 & 1 & $\frac{1}{3}= 0.3333...$ & 0.330445 \\
\hline
1000 & 0.001 & 0 & 2 & $\frac{1}{3}= 0.3333...$ & 0.321854 \\
\hline
100 & 0.001 & 1 & 1 & $\frac{3}{5}= 0.6$ & 0.600363 \\
\hline
100 & 0.001 & 2 & 1 & $\frac{5}{7}= 0.714285....$ & 0.714286 \\
\hline
100 & 0.001 & 3 & 1 & $\frac{7}{9}= 0.7777...$ & 0.777777\\
\hline
100 & 0.001 & 4 & 1 & $\frac{9}{11}= 0.8181...$ & 0.818176\\
\hline
100 & 0.3 & 49 & 1 & $\frac{99}{101}= 0.98019801...$ & 0.980189\\
\hline
\end{tabular}
\caption{Numerical results for the classical Li\'enard equation given in \eqref{Application 1 system}.}
\label{table classical}
\end{table}

\subsection{A two-stroke oscillator}\label{section-examp2}
In this section we consider a two-stroke oscillator studied in \cite{Martin}:
 \begin{equation}
\label{Application 2 system}
    \begin{vf}
        \dot{x} &=& y(\delta -y) \\
        \dot{y} &=&(-x+\alpha y)\cdot (\delta -y) - \epsilon\big(\beta - \gamma x\big),
    \end{vf}
\end{equation}
where $\alpha,\beta, \gamma, \delta >0$ and $\epsilon\ge 0$ is the singular perturbation parameter. Following \cite{Martin} or \cite{CriticalityHopf}, we deal with a slow-fast Hopf point (in a non-standard form) in \eqref{Application 2 system} at $p=(\alpha \delta, \delta)$, for $\beta=\alpha \gamma \delta$. When $\epsilon=0$, \eqref{Application 2 system} has the critical curve $\{y=\delta\}$: $p$ is a nilpotent singularity and it separates the normally attracting branch $x<\alpha\delta$ and the normally repelling branch $x>\alpha\delta$. It is not difficult to see that the slow dynamics along the critical curve is given by 
\[\frac{dx}{ds}=\frac{\beta\delta-\gamma\delta x}{\alpha\delta-x}, \ \ x\ne \alpha\delta,\]
where $s=\epsilon t$ is the slow time and $t$ is the fast time attached to \eqref{Application 2 system}. We assume that $\beta=\alpha \gamma \delta$. Then the above differential equation becomes $\frac{dx}{ds}=\gamma\delta>0$. Thus, we deal with a regular slow dynamics (including the contact point $p$) pointing from the attracting branch to the repelling branch. The main difference between \eqref{Application 1 system}, given in a standard form, and \eqref{Application 2 system} is that in \eqref{Application 2 system} the fast foliation is not given by horizontal lines. Near $p$, we have concave down parabola like fast movements and they have a quadratic contact with the critical curve at $p$ (this can be easily seen from \eqref{Application 2 system} with $\epsilon=0$). We take a section $\sigma\subset\{x=\alpha\delta\}$, parameterized by $y\sim \delta$ ($y=\delta$ corresponds to $p$).
For $y=\tilde y_0\in\sigma$ and $\tilde y_0>\delta$, we generate a sequence $(\tilde y_k)_{k\ge 0}$ using
\begin{equation}\label{Application 2 integral}
I(\tilde y_k,\tilde y_{k+1})=\frac{1}{\gamma\delta }\int_{\alpha(\tilde{y}_k)}^{\omega(\tilde{y}_{k+1})}(x-\alpha \delta) dx=0,\nonumber
\end{equation}
where $\alpha(y)>\alpha\delta$ (resp. $\omega(y)<\alpha\delta$) is the $\alpha$-limit point (resp. $\omega$-limit point) of the fast orbit of \eqref{Application 2 system} with $\epsilon=0$ passing through $y\in \sigma$, with $y>\delta$. We numerically compute $\alpha(y)$ and $\omega(y)$, for a given $y>\delta$.  
Because $\tilde{y}_k \rightarrow \delta$, we use the translation $y_k=\tilde{y}_k-\delta$, $\forall k \ge 0$, and $y_k\to 0$.   

\smallskip

Table \ref{table two-stroke} shows the estimated Minkowski dimension of $(y_k)_{k\ge 0}$ (i.e. $(\tilde y_k)_{k\ge 0}$) for different positive parameters $\alpha,\gamma,\delta$ and $\beta=\alpha \gamma \delta$.
\begin{table}[]
\centering
\begin{tabular}{|c|c|c|c|c|c|c|c|c|c|c|c|}
\hline
\textbf{\# Iterations} &  \textbf{$\tilde{y}_0$} &  \textbf{$\alpha$} & \textbf{$\delta$} & \textbf{$\gamma$} & \textbf{$\beta$} &  \textbf{Theoretical Value} & \textbf{Results} \\
\hline
1000 & 1.1 & 1 & 1 & 1 & 1	& $\frac{1}{3}= 0.3333...$ & 0.335137 \\
\hline
1000 & 1.1 & 1 & 1 & 10 & 10 & $\frac{1}{3}= 0.3333...$ & 0.335137  \\
\hline
1000 & 1.1 & 2 & 1 & 1	& 2	& $\frac{1}{3}= 0.3333...$ & 0.324280 \\
\hline
1000 & 10.1 & 5 & 10 & 1 & 50 & $\frac{1}{3}= 0.3333...$ & 0.331570  \\
\hline
\end{tabular}
\caption{Numerical results for the two-stroke oscillator.}
\label{table two-stroke}
\end{table}
The numerical results indicate that $\dim_B\mathcal (y_k)_{k\ge 0}=\frac{1}{3}$. From Definition \ref{def-intrin-codimen} it follows that the fractal codimension of $p$ is $1$ (note that $j=\frac{3\dim_B(y_k)_{k\ge 0}-1}{2(1-\dim_B(y_k)_{k\ge 0})}=0$). Following Theorem \ref{theorem-Lienard}.1 or Theorem \ref{theorem-Lienard}.3, the slow-fast Hopf point $p$ can produce at most $1$ limit cycle. This conclusion, based on our numerical results, can be theoretically proven (see \cite{CriticalityHopf,Martin}).  

\subsection{Example  with a general contact order and singularity order}\label{section-examp4}
In this section we deal with a slow-fast system of type \eqref{normal form 2}:
\begin{equation}
\label{Application 4 system}
    \begin{vf}
        \dot{x} &=& y-x^n \\
        \dot{y} &=&\epsilon \big(\beta x^m+ \alpha x^{m+2j+1} \big),
    \end{vf}
\end{equation}
 where $n\ge 2$ is even, $m\ge 1$ is odd, $j\in \mathbb{N}_0$, $\alpha \ne 0$ and $\beta= \pm 1$. We have that $m\leq 2  (n-1)$. We numerically verify some  Minkowski dimension results in Theorem \ref{theorem-1} for fractal sequences. We define a section $\sigma\subset \{x=0\}$ as in Section \ref{subsectionNintrinsic} (see Fig. \ref{fig-sequence}). 
 The slow divergence integral \eqref{SDI-contact} becomes
 \[I(y,\tilde y)=-\int_{-y^{1/n}}^{\tilde y^{1/n}}\frac{(nx^{n-1})^2}{\beta x^m+ \alpha x^{m+2j+1}}dx.\]

Again, for a $y_0>0$ small enough we generate a sequence $(y_k)_{k\ge 0}$ tending to zero using the entry-exit relation, $I(y_k,y_{k+1})=0$ or $I(y_{k+1},y_{k})=0$, depending on the sign of $\alpha$ and $\beta$ (see Section \ref{subsectionNintrinsic}).

Following Theorem \ref{theorem-1}, we know that 
    $\dim_B\mathcal (y_k)_{k\ge 0}= \frac{2j+1}{n+2j+1}$.
In Table \ref{table general} one can find numerically computed Minkowski dimensions. This corresponds to our theoretical expectations. Just as expected the parameters $\alpha, \beta$ and $m$ do not influence the Minkowski dimension. 
\begin{table}[]
\centering
\begin{tabular}{|c|c|c|c|c|c|c|c|c|c|c|c|c|}
\hline
\textbf{\# It} &  \textbf{${y}_0$} &  \textbf{m} & \textbf{n} & \textbf{j}& \textbf{$\alpha$} & \textbf{$\beta$} &  \textbf{Theoretical Value} & \textbf{Results} \\
\hline
2000 & 0.1 & 1 & 2 & 0 & 1 & 1 & $\frac{1}{3}=0.3333...$ &  0.345550 \\
\hline
2000 & 0.1 & 1 & 2 & 0 & 1 & -1 & $\frac{1}{3}=0.3333...$ & 0.345550 \\
\hline
2000 & 0.1 & 1 & 2 & 10 & 1 & 1 & $ \frac{21}{23}= 0.913043...$ & 0.920386 \\
\hline
2000 & 0.1 & 1 & 4 & 10	& 1	& 1	& $\frac{21}{25}=0.84$ & 0.858920 \\
\hline
2000 & 0.1 & 1 & 10 & 10 & 1 & 1 & $\frac{21}{31}= 0.677419...$ & 0.673676\\
\hline
2000 & 0.1 & 3 & 4 & 10 & 1 & 1 & $\frac{21}{25}=0.84$ & 0.858265 \\
\hline
2000 & 0.1 & 9 & 10 & 5	& 5	& 1 & $\frac{11}{21}= 0.523809...$ & 0.523656 \\
\hline
2000 & 0.1 & 99 & 100 & 50 & 1 & 1 & $\frac{101}{201}= 0.502487...$ & 0.502158 \\
\hline
\end{tabular}
\caption{Numerical results for \eqref{Application 4 system} with general contact order and singularity order. (It=Iterations.)}
\label{table general}
\end{table}

\section*{Acknowledgments}

This research was supported by: Croatian Science Foundation (HRZZ) grant PZS-2019-02-3055 from “Research Cooperability” program funded by the European Social Fund.
Additionally, the research of Goran Radunovi\'c was partially supported by the HRZZ grant UIP-2017-05-1020.

\bibliographystyle{plain}
\bibliography{bibtex}
\end{document}